\documentclass[letterpaper]{article} 
\usepackage{aaai23}  
\usepackage{times}  
\usepackage{helvet}  
\usepackage{courier}  
\usepackage[hyphens]{url}  
\usepackage{graphicx} 
\usepackage{natbib}  
\usepackage{caption} 
\usepackage{algorithm}
\usepackage{algorithmic}
\usepackage{amsmath}
\usepackage{indentfirst}
\usepackage{amsfonts}
\usepackage{multirow}
\usepackage{stmaryrd}

\newtheorem{theorem}{Theorem}
\newtheorem{lemma}{Lemma}

\newtheorem{remark}{Remark}

\usepackage{newfloat}
\usepackage{listings}

\usepackage{booktabs}
\usepackage{subfigure} 
\usepackage{amssymb}
\usepackage{mathrsfs}
\usepackage{cite}

\DeclareCaptionStyle{ruled}{labelfont=normalfont,labelsep=colon,strut=off} 
\floatstyle{ruled}
\newfloat{listing}{tb}{lst}{}
\floatname{listing}{Listing}
\title{Stepdown SLOPE for Controlled Feature Selection}
\author{
    Jingxuan Liang\textsuperscript{\rm 1}, Xuelin Zhang \textsuperscript{\rm 2}, Hong Chen\textsuperscript{\rm 1, 4, 6,}\thanks{Corresponding author.}, Weifu Li\textsuperscript{\rm 1, 5, 6}, Xin Tang\textsuperscript{\rm 3}\\
}
\affiliations{
    \textsuperscript{\rm 1}College of Science, Huazhong Agricultural University, Wuhan 430070, China\\
    \textsuperscript{\rm 2} College of Informatics, Huazhong Agricultural University, Wuhan 430070, China\\
      \textsuperscript{\rm 3}Ping An Property \& Casualty Insurance Company, Shenzhen, China\\
    \textsuperscript{\rm 4}Engineering Research Center of Intelligent Technology for Agriculture, Ministry of Education, Wuhan 430070, China\\
    \textsuperscript{\rm 5}Key Laboratory of Smart Farming for Agricultural Animals, Wuhan 430070, China\\
       \textsuperscript{\rm 6}Hubei Engineering Technology Research Center of Agricultural Big Data,  Wuhan 430070, China\\
     chenh@mail.hzau.edu.cn

}

\usepackage{bibentry}

\begin{document}
\maketitle

\begin{abstract}
Sorted L-One Penalized Estimation (SLOPE) has shown the nice theoretical property as well as empirical behavior recently on the false discovery rate (FDR) control of high-dimensional feature selection by adaptively imposing the non-increasing sequence of tuning parameters on the sorted $\ell_1$ penalties. This paper goes beyond the previous concern limited to the FDR control by considering the stepdown-based SLOPE to control the probability of $k$ or more false rejections ($k$-FWER) and the false discovery proportion (FDP). Two new SLOPEs, called $k$-SLOPE and F-SLOPE, are proposed to realize $k$-FWER and FDP control respectively, where the stepdown procedure is injected into the SLOPE scheme. For the proposed stepdown SLOPEs, we establish their theoretical guarantees on controlling $k$-FWER and FDP under the orthogonal design setting, and also provide an intuitive guideline for the choice of regularization parameter sequence in much general setting. Empirical evaluations on simulated data validate the effectiveness of our approaches on controlled feature selection and support our theoretical findings. 
\end{abstract}

\section{Introduction} \label{sec1}
Feature selection aims to find the informative features from high-dimensional empirical observations, which is one of key research fields of machine learning. Typical feature selection methods include sparse linear models (e.g., Lasso \cite{lasso}), sparse additive models (e.g., SpAM \cite{spam}, GroupSAM \cite{chenh-nips17}, SpMAM \cite{chenh-tnnls2020}), tree-based models (e.g., random forest \cite{randomforest}), and sparse neural networks (e.g., LassoNet \cite{lemhadri-2021lassonet}). 

Following this line, the controlled feature selection further addresses the selection quality with low false discovery rate (FDR) guarantee, which has attracted the increasing attention recently due to its wide applications, e.g.,  in bioinformatics and biomedical \cite{aggarwal-2016biofdr,yu-2021biofdr}. There are mainly three branches of learning systems for controlled feature selection: the multiple hypothesis test \cite{benjamini-1995bh,ferreira-2006bh,kFWER2,stepup}, the knockoffs filter \cite{barber-2015knockoffs,Cands2016PanningFG,barber-2020knockoffs,romano-2020knockoffs}, and the Sorted L-One Penalized Estimation (SLOPE) \cite{slope,slopestatics,Groupslope}. As a classic strategy for feature selection, the \citeauthor{benjamini-1995bh} (BH) procedure is formulated by jointly considering p-values of multiple hypothesis testing. Despite this procedure enjoys nice theoretical properties on the FDR control, it may face the computation challenge for nonlinear and complex regression estimation \cite{javanmard2019false}. As a novel feature filter scheme, the knockoffs inference has solid theoretical foundations and shows the competitive performance in real-word applications \cite{barber-2015knockoffs,barber-2020knockoffs,zhao-2022error,yu-2021biofdr}. Particularly, an error-based knockoffs inference framework is formulated in \cite{zhao-2022error} to further realize the controlled feature selection from the perspectives of the probability of $k$ or more false rejections ($k$-FWER) and the false discovery proportion (FDP). Different from screening out the active feature with the help of knockoff features, SLOPE focuses on the regularization design for sparse feature selection, which adaptively imposes a non-increasing sequence of tuning parameters on the sorted $\ell_1$ penalties \cite{slope,Groupslope,abslope}. 
 
Although rapid progresses on its optimization algorithm \cite{slope,Groupslope} and theoretical properties \cite{slopestatics}, all the existing works of SLOPE are limited to the FDR control only. Naturally, it is important to explore new SLOPE for controlled feature selection under other statistical criterion, e.g., $k$-FWER and FDP.
To fill this gap, we propose new SLOPE approaches, called $k$-SLOPE and F-SLOPE, to realize feature selection with the $k$-FWER and FDP control respectively. Different from the previous method relying on BH procedure, the proposed  SLOPEs depend on the stepdown procedure \cite{kFWER2}, which enjoy much feasibility and adpativity \cite{slope,slopestatics}. The main contributions of this paper are summarized as below:
\begin{itemize}
\item \emph{New SLOPEs for the kFWER and the FDP control}. We integrate the SLOPE \cite{slope} and the stepdown procedure \cite{kFWER2} into a coherent way for the $k$-FWER and FDP control and formulate the respective convex optimization problem. Similarly with the flexible knockoffs inference in \cite{zhao-2022error}, our approaches also can avoid the complex p-value calculation and can be implemented feasibility.
\item \emph{Theoretical guarantees and empirical effectiveness}. 
Under the orthogonal design setting, the $k$-FWER and FDP can be provably controlled at a prespecified level for the proposed $k$-SLOPE and $F$-SLOPE, respectively.   In addition, we provide an intuitive theoretical analysis for the choice of the regularizing sequence in general setting. Simulated experiments validate the effectiveness of our SLOPEs on the $k$-FWER and FDP control, and verify our theoretical findings. 
\end{itemize}
\section{Related Work}\label{sec2}
To better highlight the novelty of the proposed method, we review the related SLOPE methods as well as the relationship among FDR, $k$-FWER and FDP.

{\bf{SLOPE} Methods.}
SLOPE \cite{slope} can be considered as a natural extension of Lasso \cite{lasso}, where the regression coefficients are penalized according to their rank. 
One notable choice of the regularization sequence $\{\lambda_i\}$ is given by the BH \cite{benjamini-1995bh} critical values $\lambda_{\mathrm{BH}}(i)=\Phi^{-1}(1-\frac{iq}{2m})$, where $q \in (0,1)$ is the desired FDR level, $m$ is the characteristic number and $\Phi(\cdot)$ is the cumulative distribution function of a standard normal distribution. The main motivation behind SLOPE is to provide finite sample guarantees on regression estimation and FDR control,  where FDR is defined as the expected proportion of irrelevant regressors among all selected predictors. When $X$ is an orthogonal matrix, SLOPE with $\lambda_{\mathrm{BH}}$ controls FDR at the desired level in theory. Besides, a remarkable feature is that SLOPE does not require any knowledge of the degree of sparsity, yet automatically yields optimal total squared errors over a wide range of $\ell_0$-sparsity classes. 

To improve computing efficiency, a sparse semismooth Newton-based augmented Lagrangian technique was proposed to solve the more general SLOPE model \cite{luo2019solving}. A heuristic screening rule for SLOPE based on the strong rule for the lasso was first presented in order to improve the numerical procedures efficiency of SLOPE, especially in the setting of estimating a complete regularization path \cite{strongslope}. And \citet{larsson} also proposed a new fast algorithm to solve the SLOPE optimization problem, which combined proximal gradient descent and proximal coordinate descent steps. Besides the above works on algorithm optimization, there are extensive studies on SLOPE with properties \cite{slopestatics,bellec,kos}, model improvements \cite{Groupslope,lee,riccobello,abslope} and applications \cite{brzyski,kremer}. As we know, there is no any touch to address the SLOPE-based feature selection with $k$-FWER or FDP control guarantees. 

{\bf{Statistical Metrics: FDR, $k$-FWER and FDP}.}
\citet{benjamini-1995bh} formulated the BH procedure to the control the expectations of FDP, called FDR control. Then, \citet{kFWER2} proposed both the single step procedure and the stepdown procedure in order to ensure the $k$-FWER control. \citet{kFWER2} also considered the FDP control and provided  two stepdown procedures for controlling the FDP under mild conditions with the p-values dependence structure or no any dependence supposition. With the help of stepdown procedures \cite{kFWER2}, there are studies on feature selection with the $k$-FWER control \cite{stepup,romano2007control,aleman2017effects,zhao-2022error} and the FDP control \cite{stepup,romano2007control,fan2010selective,delattre2015new,zhao-2022error}. However, most of these procedures may depend on the p-values to assess the importance of each feature or the assumption of structures. Moreover, the traditional calculation of p-value relies on the large-sample asymptotic theory usually, which may  no longer be true in the setting of high-dimensional finite samples \cite{Cands2016PanningFG,Fan}. 

It is necessary to explain the relationship between FDR, FDP and $k$-FWER. Given $\gamma,\alpha\in (0,1)$, the FDP control means the $\rm{Prob}(\mathrm{FDP}>\gamma)$ at the level $\alpha$. Recall that the FDP concerns \begin{equation}
	\rm{Prob}\{\mathrm{FDP}>\gamma\}<\alpha,
\end{equation}
and FDR is the expectation of FDP, i.e., $\mathrm{FDR}=\mathbb{E}\mathrm{(FDP)}$. It is easy to verify that 
\begin{small}
\begin{equation*}
\begin{aligned}
\mathrm{FDR}
 \leq \gamma \rm{Prob}\{\mathrm{FDP} \leq \gamma\}+\rm{Prob}\{\mathrm{FDP}>\gamma\},
\end{aligned}
\end{equation*}
\end{small}
and then
\begin{equation*}
\frac{\mathrm{FDR}-\gamma}{1-\gamma} \leq \rm{Prob}\{\mathrm{FDP}>\gamma\} \leq \frac{\mathrm{FDR}}{\gamma},
\end{equation*}
where the last inequality follows from the Markov’s inequality. Clearly, if a method controls FDR at level $q$, then it also controls $\mathrm{FDP}\leq q/\gamma$. Conversely, if the FDP is controlled, i.e. $\rm{Prob}(\mathrm{FDP}>\gamma)<\alpha$, and then the FDR is bounded by $(1-\gamma)\alpha+\gamma$. Therefore, a procedure with the FDP control often can control the FDR \cite{vanderLaanDudoitPollard}. Furthermore, \citet{Farcomen} pointed out that, compared with the FDR control, the $k$-FWER control is more desirable when powerful selection results can be made.

\section{Preliminaries}
\label{sec3}
This section recalls some necessary backgrounds involved in this paper, e.g., SLOPE \cite{slope} and the stepdown procedure \cite{kFWER2}. The main notations used in this paper are summarized in \emph{Supplementary Material A}.  

\subsection{Problem Formulation}
Let $\mathcal{X}\subset\mathbb{R}^m$ and $\mathcal{Y}\subset \mathbb{R}$ be the compact input space and corresponding output space, respectively. Consider samples $\{(x_i,y_i)\}_{i=1}^n$ independently drawn from an unknown distribution on $\mathcal{X} \times \mathcal{Y}$. Denote $$X:=(X_1,X_2,\cdots,X_n)^T\subset\mathbb{R}^{n\times m}$$
and $$y=(y_1,y_2,\cdots,y_n)^T\in\mathbb{R}^n.$$
The module length of each column vector of $X$ is equal to 1. The output vector $y$ is generated by the following multiple linear regression model:
\begin{align}
	y = X\beta+\epsilon,\label{linear}
\end{align}
where $\beta\in\mathbb{R}^m$ represents the coefficient vector and $\epsilon \sim N(0,\sigma^2I_n)$. In sparse high-dimensional  regression, we often assume that $\beta$ satisfies a sparse structure. Let $V$ be the number of false selected features and let $R$ be total number of identified features. The FDP, FDR and $k$-FWER are respectively defined as 

\begin{align*}
    &\mathrm{FDP}=\frac{V}{\max\{R,1\}},~~ \mathrm{FDR}=\mathbb{E}(\mathrm{FDP})
\end{align*}
and
\begin{align*}
	&k\text{-}\mathrm{FWER}={\rm{Prob}}\{V\geq k\}. 
\end{align*}	

\subsection{SLOPE}
SLOPE is proposed by \citet{slope} for controlled feature selection in high dimensional sparse cases, which replaces the $\ell_1$ penalty in Lasso \cite{lasso} with the sorted $\ell_1$ penalty. 
The learning scheme of SLOPE \cite{slope} is formulated as 
\begin{align}
	\mathop{\arg\min}_{\beta\in \mathbb{R}^m}\frac{1}{2}||y-X\beta||^2+\sum_{i=1}^m\lambda_i|\beta|_{(i)},\label{SLOPEeq}
\end{align}
where 
the regularization parameters $\lambda_{1}\geq\cdots\geq\lambda_{m}\geq0$ and the regression coefficients $|\beta|_{(1)}\geq |\beta|_{(2)}\geq\cdots\geq |\beta|_{(m)}$ are all non-negative non-decreasing sequences. When $\lambda_1=\lambda_2=\cdots=\lambda_m$, the optimizing scheme (\ref{SLOPEeq}) obviously reduces to the Lasso \cite{lasso}. Given a desired level $q$, SLOPE controls FDR using the sequence of parameters $\mathrm{\lambda_{BH}=\{\lambda_{BH}(1),\lambda_{BH}(2),\cdots,\lambda_{BH}}(m)\}$ with 
\begin{equation}
	\mathrm{\lambda_{BH}}(i)=\sigma\cdot\Phi^{-1}(1-\frac{i q}{2m})\label{lambdaBH},
\end{equation}
where $\Phi(\cdot)$ denotes the cumulative distribution function of the standard normal distribution under orthogonal design. 
\begin{theorem}\cite{slope}
In the linear model with the orthogonal design $X$ and $\epsilon\sim N (0,\sigma^2 I_n)$, the SLOPE \eqref{SLOPEeq} with the regularization parameter sequence  (\ref{lambdaBH}) satisfies 
$$
\mathrm{FDR} \leq q \frac{m_{0}}{m},
$$
where $m_0$ is the number of true null hypotheses and $q$ is the desired FDR level.
\end{theorem}

Theorem 1 illustrates the theoretical guarantee of FDR control for SLOPE equipped with $\lambda_{\mathrm{BH}}$ induced by the BH procedure \cite{benjamini-1995bh}. In this paper, we are not limited to the FDR control, but extend to the $k$-FWER and FDP control by replacing the BH procedure with the stepdown procedure \cite{kFWER2}.
\begin{algorithm}[ht]
\caption{Accelerated proximal gradient algorithm for SLOPE (\ref{SLOPEeq})}
\textbf{Input}: Training set $X\in\mathbb{R}^{n\times m}$ and $y\in\mathbb{R}^n$ and parameter $\lambda=(\lambda_1,\lambda_2,...,\lambda_m)$.\\
\textbf{Initialization}: $a^0\in\mathbb{R}^m$, $b^0=a^0$ and $\theta_0=1$.
\begin{algorithmic}[0] 
\FOR {k = 0,1,$\cdots$}
\STATE $b^{k+1}=\operatorname{prox}_{t_{k} J_{\lambda}}\left(a^{k}-t_{k} X^{\prime}\left(X a^{k}-y\right)\right)$
\STATE $\theta_{k+1}^{-1}=\frac{1}{2}\left(1+\sqrt{\left.1+4 / \theta_{k}^{2}\right)}\right.$
\STATE $a^{k+1}=b^{k+1}+\theta_{k+1}\left(\theta_{k}^{-1}-1\right)\left(b^{k+1}-b^{k}\right)$
\ENDFOR
\end{algorithmic}
\textbf{Output}: $a$ satisfying the stopping criteria.
\label{accslope}
\end{algorithm}

From the computing side, the optimization objective function of SLOPE (\ref{SLOPEeq}) is convex but non-smooth, which can be implemented efficiently by the proximal gradient descent algorithm \cite{slope}. For completeness, we state the computing steps of SLOPE in Algorithm \ref{accslope}, which also suits for our variants of SLOPE. Here, $J_\lambda=\sum_{i=1}^m\lambda_i|\beta|_{(i)}$ and the step lengths get by backtracking line search and satisfy $t_k<2/||X||^2$ \cite{beck2009fast,becker2011templates}. Moreover, \citet{slope} also derive concrete stopping criteria through duality theory.

\subsection{Stepdown Procedure}
The stepdown procedure \cite{kFWER2} aims to control $k$-FWER and FDP, i.e., given $\alpha,r\in (0,1)$, 
\begin{align}
	&k\text{-}\mathrm{FWER}\leq\alpha
	\label{kFWERcantrol}
\end{align}
and
\begin{align}
	&\rm{Prob}\{\mathrm{FDP}>\gamma\}\leq\alpha.
	\label{FDPcontrol}
\end{align}
Suppose that there are $m$ individual tests $H_{1},...,H_m$, whose corresponding p-values are $\hat{p}_1,...,\hat{p}_m$. Let $\hat{p}_{(1)}\leq \hat{p}_{(2)}\leq...\leq \hat{p}_{(m)}$ be the ordered p-values and let the non-negative non-decreasing sequence $\alpha_1\leq\alpha_2...\leq \alpha_m$ be the $k$-FWER thresholds. The hypotheses corresponding to the sorted p-values are defined as $H_{(1)},H_{(2)}...,H_{(m)}$. Then the stepdown procedure is defined stepwise as follows:

$Step\,0$: Let $i=0$.

$Step\,1$: If $\hat{p}_{(i+1)}\geq \alpha_{i+1}$, go to step 2. Otherwise, set $i=i+1$ and repeat $Step$ 1.

$Step\,2$: Reject $H_{(j)}$ for $j\leq k$ and accept $H_{(j)}$ for $j>k$.
In other words, if $p_{(1)}>\alpha_1$, no null hypotheses are rejected. Otherwise, if $H_{(1)},H_{(2)}...,H_{(r)}$ are rejected, the largest $r$ satisfies
\begin{equation}
	p_{(1)}\leq \alpha_1,p_{(2)}\leq\alpha_2,...,p_{(r)}\leq\alpha_r.\label{down}
\end{equation} 
Based on the stepdown procedure, \citet{kFWER2} provided two different thresholds to ensure the $k$-FWER control and the FDP control, respectively.
\begin{theorem}\cite{kFWER2}
For testing $H_i,i= 1,...,m$, given $k$ and $\alpha\in(0,1)$, the stepdown procedure described in (\ref{down}) with \begin{equation}
\alpha_{i}= \begin{cases}\frac{k \alpha}{m}, & i \leq k \\ \frac{k \alpha}{m+k-i}, & i>k\end{cases}	\label{kFWER_threlod}
\end{equation}
controls the $k$-FWER, that is, (\ref{kFWERcantrol}) holds.

\end{theorem}

\begin{theorem}\cite{kFWER2}
For testing $H_i,i= 1,...,m$, given $\alpha,\gamma\in (0,1)$, if the p-values of false null hypotheses are independent of the true ones, the stepdown procedure described in (\ref{down}) with 
\begin{equation}
\alpha_{i}=\frac{(\lfloor\gamma i\rfloor+1) \alpha}{m+\lfloor\gamma i\rfloor+1-i}\label{FDP_threlod}
\end{equation}
controls the FDP in the sense of (\ref{FDPcontrol}).
\end{theorem}

Theorems 2 and 3 demonstrate that the stepdown procedure enjoys the theoretical guarantees on the $k$-FWER control and FDP control under ingenious selections of $\alpha_i$.  Indeed, these theoretical properties of stepdown procedure motivate our designs for new SLOPE algorithms. 

\section{Methodology}\label{sec4}
This section injects the stepdown procedure \cite{kFWER2} into the classical SLOPE \cite{slope} to formulate our new stepdown SLOPEs for controlled feature selection to ensure the $k$-FWER control and the FDP control. Here, we provide the sequences of tuning parameters under the orthogonal design for the $k$-FWER control and the FDP control, respectively. Furthermore, we present an intuitive theoretical analysis for the selection of regularization parameters in general setting.

\subsection{Orthogonal Design}
It has been illustrated in \citet{slope} that the linking between multiple tests and model selection for SLOPE under the orthogonal design. Following this line, we assume that $X$ is an $n\times m$ dimensional orthogonal matrix, i.e, $X'X=I_m$ and $\epsilon \sim N(0,\sigma^2I_n)$ is an $n$-dimensional column vector with known variance. Then, the linear regression model $$y=X\beta+\epsilon$$
is transformed into
\begin{equation*}
	\tilde{y}=X'y=\beta+X'\epsilon\sim N(\beta,\sigma^2I_p).
\end{equation*}
It is well known that the problem of selecting effective features can be  simplified as a multiple hypothesis test problem. Denote $m$ hypotheses as $H_i:\beta_i=0,1\leq i\leq m$. If $H_i$ is rejected, $\beta_i$ is considered as an effective feature and vice versa. \citet{slope} gave the selection mechanism of regularization parameters for SLOPE through the BH procedure \cite{benjamini-1995bh} under the orthogonal design. For brevity, we call the  proposed methods as $k$-SLOPE and F-SLOPE with respect to the control of $k$-FWER and FDP, respectively. 

The regularization scheme of $k$-SLOPE is formulated as
\begin{equation}
\mathop{\arg \min}_{\beta\in \mathbb{R}^m}\frac{1}{2}\|y-X \beta\|_{l_2}^{2}+ \sigma\cdot\sum_{i=1}^m\lambda_{k\text{-}\mathrm{FWER}}(i)|\beta|_{(i)},\label{kFWER_model}
\end{equation} 
where
\begin{equation} \label{sequence1}
\lambda_{k\text{-}\mathrm{FWER}}(i)= \begin{cases}\Phi^{-1}(1-k \alpha/{2m}), & i \leq k \\ \Phi^{-1}(1-k\alpha/{2(m+k-i)}), & i>k.\end{cases}
\end{equation}
The $k$-SLOPE equipped with \eqref{sequence1}  yields the following theoretical property, which has been proved in \emph{Supplementary Material B}. 
\begin{theorem}
	In the linear model (\ref{linear}) with the orthogonal matrix $X$ and noise $\epsilon\sim N(0,\sigma^2I_n)$, given $k$ and $\alpha\in(0,1)$, the $k$-FWER of the $k$-SLOPE model (\ref{kFWER_model}) satisfies (\ref{kFWERcantrol}).
\end{theorem}
Theorem 4 illustrates that $k$-SLOPE controls the $k$-FWER under the orthogonal design, which has been proved in \emph{Appendix}. Although the $\lambda_{k\text{-}\mathrm{FWER}}(i)$'s are chosen with reference \cite{kFWER2}, (\ref{kFWER_model}) is not equivalent to the stepdown procedure described above. We also empirically support this theoretical guarantee by experimental analysis.

Generally, the number of false selected features that people are willing to abide is directly proportional to the number of identified features. Therefore, we may be no longer concerned about $k$-FWER, but about FDP. Similar to (\ref{kFWER_model}), the convex optimization problem of F-SLOPE is formulated as
\begin{equation}
\mathop{\arg \min}_{\beta\in \mathbb{R}^m}\frac{1}{2}\|y-X \beta\|_{l_2}^{2}+ \sigma\cdot\sum_{i=1}^m\lambda_{\mathrm{FDP}}(i)|\beta|_{(i)},\label{FDP_model}
\end{equation}
where
\begin{align*}
&\lambda_{\mathrm{FDP}}(i)=\Phi^{-1}(1-\frac{(\lfloor\gamma i\rfloor+1) \alpha}{2(m+\lfloor\gamma i\rfloor+1-i)}).
\end{align*}
The selection of regularization parameters also produces the following theoretical guarantee. 

\begin{theorem}\label{fdp_t}
	In the linear model (\ref{linear}) with the orthogonal matrix $X$ and noise $\epsilon\sim N(0,\sigma^2I_n)$, given $\alpha,\gamma\in (0,1)$, the FDP of the F-SLOPE model (\ref{FDP_model}) satisfies (\ref{FDPcontrol}).
\end{theorem}
Theorem 5 assures the ability of FDP control for F-SLOPE under the orthogonal design setting, which has been established in \emph{Supplementary Material B}. The only difference between the F-SLOPE model (\ref{FDP_model}) and the $k$-SLOPE model (\ref{kFWER_model}) is the selection mechanism of the sequence for penalty parameters. The conclusion is also supported by the later orthogonal experiments. Moreover, the optimization algorithm of $k$-SLOPE and F-SLOPE is the same as that of SLOPE because they are all convex and non-smooth. More optimization details are present in Algorithm \ref{accslope}.

\begin{table*}[ht]
\centering
\renewcommand\arraystretch{1.1}{
\setlength{\tabcolsep}{2.5mm}{
\begin{tabular}{c|ccc|ccc|ccc}
\hline
\multirow{2}{*}{$t$} & \multicolumn{3}{c|}{SLOPE}            & \multicolumn{3}{c|}{$k$-SLOPE}          & \multicolumn{3}{c}{F-SLOPE}           \\
                   & $\mathrm{Prob(FDP}>\gamma)$ & FDR   & Power & $\mathrm{Prob(FDP}>\gamma)$ & FDR   & Power & $\mathrm{Prob(FDP}>\gamma)$ & FDR   & Power \\ \hline
50                 & 0.450                 & 0.094 & 1.000 & 0.001                 & 0.006 & 1.000 & 0.003                 & 0.007 & 1.000 \\
100                & 0.330                 & 0.092 & 0.995 & 0.000                 & 0.002 & 0.998 & 0.002                 & 0.005 & 1.000 \\
200                & 0.140                 & 0.080 & 0.999 & 0.001                 & 0.001 & 1.000 & 0.000                 & 0.007 & 1.000 \\
300                & 0.000                 & 0.070 & 1.000 & 0.002                 & 0.001 & 1.000 & 0.000                 & 0.005 & 0.995 \\
400                & 0.000                 & 0.058 & 1.000 & 0.000                 & 0.001 & 0.995 & 0.001                 & 0.004 & 0.994 \\
500                & 0.000                 & 0.050 & 1.000 & 0.000                 & 0.001 & 0.997 & 0.000                 & 0.005 & 0.997 \\ \hline
\end{tabular}}}
\caption{Results for controlled feature selection under the orthogonal design (different $t$ and fixed $k=5$). }
\label{tab1}
\end{table*}

\subsection{General Setting}
Usually, SLOPE is difficult to establish  solid theoretical guarantees for the FDR control in non-orthogonal setting \cite{slope}. Hence, $k$-SLOPE and F-SLOPE may also face the degraded performance under such general setting. Fortunately, \citet{slope} used their own qualitative insights to make an intuitive adjustment to the regularization parameter sequence and showed the empirical effectiveness. Analogous to SLOPE, we give the regularization parameter forms of $k$-SLOPE and F-SLOPE through theoretical analysis in general setting.

Assume $k$-SLOPE and F-SLOPE correctly detect these features and correctly estimate the signs of the regression coefficients. Let $X_S$ and $\beta_S$ be the subset of variables associated to $\beta_i\neq 0$ and the value of their coefficients, respectively. The nonzero components estimator is approximated by
\begin{equation}
	\hat{\beta}_S\approx(X_S'X_S)^{-1}(X_S'y-\lambda_S)=\hat{\beta}_{\mathrm{LSE}}-(X_S'X_S)^{-1}\lambda_S,\label{1}
\end{equation}
where $\lambda_S=(\lambda_1,...,\lambda_{|S|})'$ and $\hat{\beta}_{\mathrm{LSE}}$ is the least-squares estimator of $\beta_S$. Inspired by \cite{slope}, we calculate the distribution of $X_i'X_S(\beta_S-\hat{\beta_S})$ to determine the specific forms of the regularization parameters for $k$-SLOPE and F-SLOPE. In light of (\ref{1}), $$\mathbb{E}(\beta_S-\hat{\beta}_S)\approx (X_S'X_S)^{-1}\lambda_S$$
and
\begin{equation*}
	\mathbb{E}X_i'X_S(\beta_S-\hat{\beta}_S)\approx\mathbb{E}X_i'X_S(X_S'X_S)^{-1}\lambda_S.
\end{equation*}
Under the gaussian design, where each element of $X$ is i.i.d $N(0,1/n)$, 
\begin{align*}
	\mathbb{E}(X_i'X_S(X'_SX_S)^{-1}\lambda_S)^2&=\frac{1}{n}\lambda_S'\mathbb{E}(X_S'X_S)^{-1}\lambda_S\\
	&=w(|S|)\cdot||\lambda_S||^2,\notag
\end{align*}
and
\begin{align*}
    w(|S|)&=\frac{1}{n-|S|-1},
\end{align*}
where $|S|$ is the number of elements of $S$,  $i\notin S$ and the second equation relies on the fact that the expected of an inverse $|S|\times|S|$ Wishart matrix with $n$ degrees of freedom is equal to $I_{|S|}/(n-|S|-1)$ \cite{wishart}. 

The $k$-SLOPE begins with $\lambda_{k\mathrm{G}}=\lambda_{k\text{-}\mathrm{FWER}}(1)$. Then, we take into account the slight increase in variance so that 
\begin{equation*}
	\lambda_{k\mathrm{G}}(2)=\lambda_{k\text{-}\mathrm{FWER}}(2)\sqrt{1+w(2)\lambda_{k\mathrm{G}}(1)^2}.
\end{equation*}
Thus, the sequence of $\lambda_{k\mathrm{G}}$ can be expressed as
\begin{equation}
	\lambda_{k\mathrm{G}}(i)=\lambda_{k\text{-}\mathrm{FWER}}(i)\sqrt{1+w(i-1)\sum_{j< i}\lambda_{k\mathrm{G}}(i)^2}.\label{kgeneral}
\end{equation}
The only difference between F-SLOPE and the $k$-SLOPE is the selection of the coefficient sequence of the penalty term. Similar with (\ref{kgeneral}), F-SLOPE starts with $\lambda_{\mathrm{FG}}=\lambda_{\mathrm{FDP}}(1)$, and then
\begin{equation}
	\lambda_{\mathrm{FG}}(i)=\lambda_{\mathrm{FDP}}(i)\sqrt{1+w(i-1)\sum_{j< i}\lambda_{\mathrm{FG}}(i)^2}.\label{fgeneral}
\end{equation}

\begin{figure}[ht]
	\centering
	\includegraphics[width=1\linewidth]{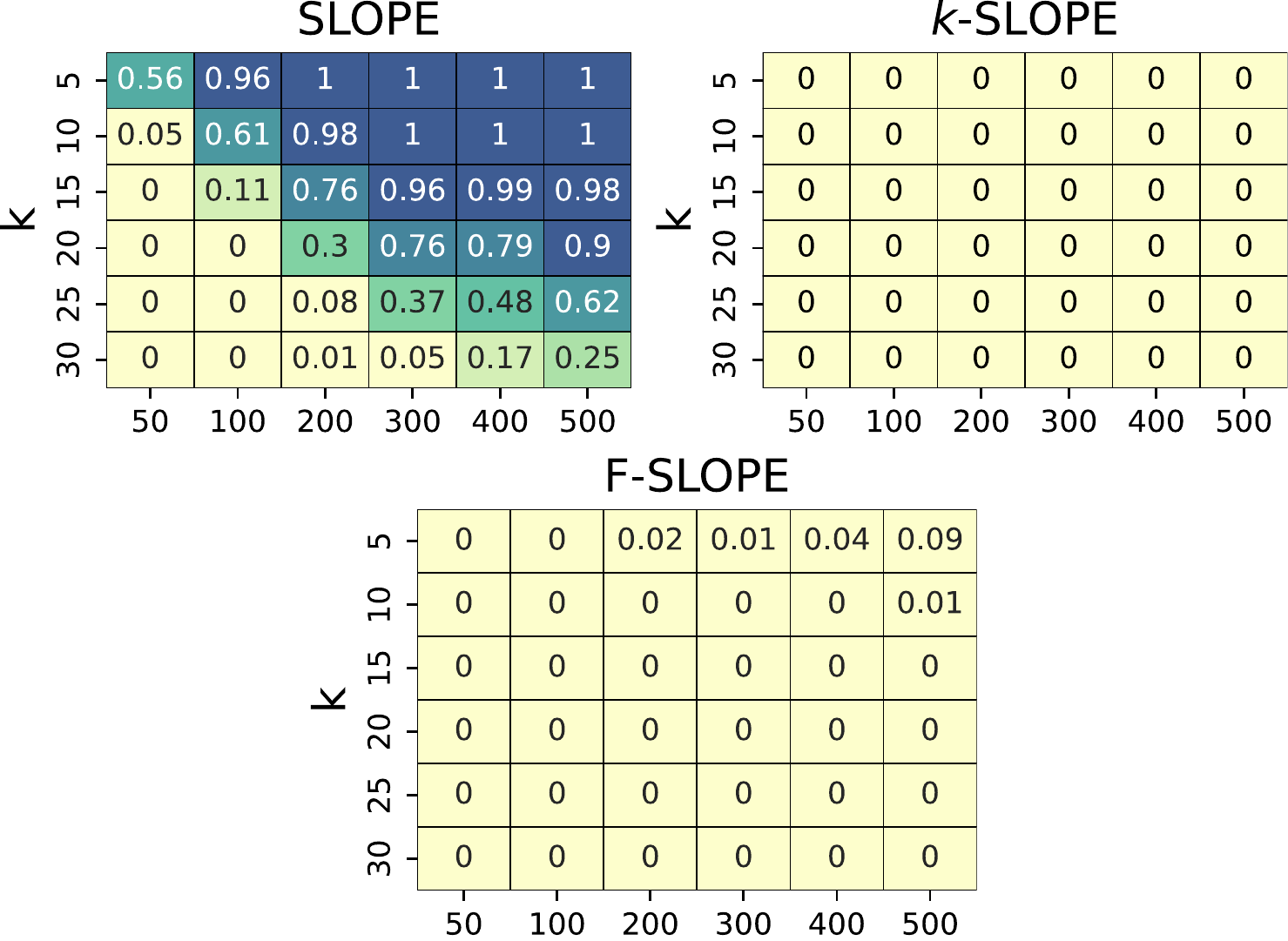}
	\caption{$k$-FWER provided by different approaches for controlled feature selection under orthogonal design (with different $k$ and $t$). The value in the small square is the size of $k$-FWER. The darker the color, the larger the $k$-FWER and vice versa.}
	\label{orth_kfwer_figure}
\end{figure}
If the coefficient sequence of the penalty term is an incremental sequence, $k$-SLOPE and F-SLOPE no longer are the convex optimization problems. Denote $k^*:=k(n,m,\alpha)$ as the subscript of global minimum, $k$-SLOPE and F-SLOPE respectively work with
\begin{equation}
\lambda_{\mathrm{kG}^{\star}}(i)= \begin{cases}\lambda_{\mathrm{kG}}(i), & i \leq k^{\star}, \\ \lambda_{k G}\left(k^{\star}\right), & i>k^{\star},\end{cases}
\end{equation}
 with $\lambda_{k\mathrm{G}} (i)$ given in (\ref{kgeneral}) and 
\begin{equation}
\lambda_{\mathrm{FG}}(i)= \begin{cases}\lambda_{\mathrm{FG}}(i), & i \leq k^{\star}, \\ \lambda_{F G}\left(k^{\star}\right), & i>k^{\star},\end{cases}
\end{equation}
with $\lambda_{\mathrm{FG}}(i)$ defined in (\ref{fgeneral})). 
When the design matrix isn't Gaussian or that columns aren't independent, we can employ the Monte Carlo estimate of the correction \cite{hammersley1954poor} instead of $w(i-1)\sum_{j<i}\lambda(i)^2$ in the formulas (\ref{kgeneral}) and (\ref{fgeneral}). 

\section{Empirical Validation}\label{sec5}
All experiments are implemented in Python on a Macbook Pro with Apple M1 and 16 GB memory. 
The reported results are the average values after repeating 100 times for each experiment.

\subsection{Experiments of Orthogonal Design Setting}

\begin{figure*}[htbp]
	\centering
	\includegraphics[width=1\linewidth]{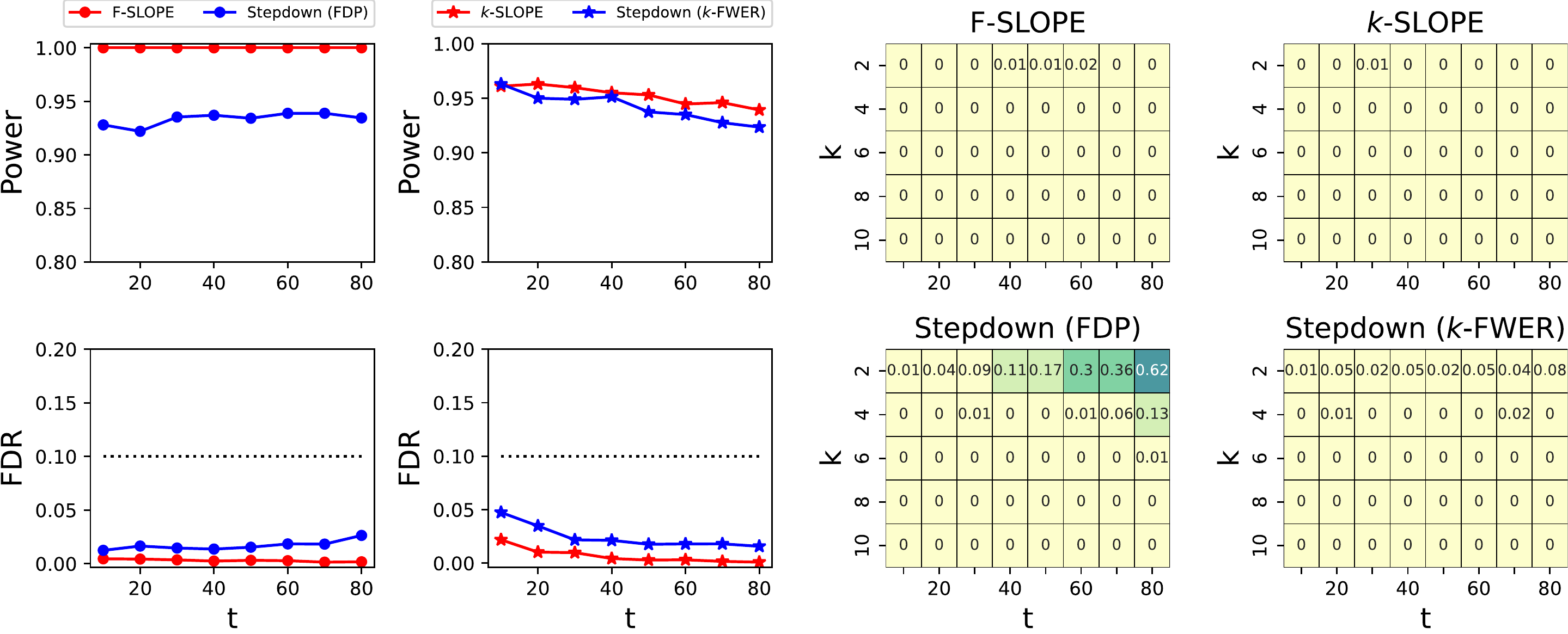}
	\caption{Result for controlled feature selection on the simulated data. The black dashed lines indicate the target FDR level. Constance for $k$-SLOPE is $k=6$ in the second column (from left to right). The value in the small square is the size of $k$-FWER in the third and fourth columns (from left to right). The darker the color, the larger the $k$-FWER and vice versa. }
\label{figure2}
\end{figure*}

Inspired by \cite{slope,Groupslope}, we draw the design matrix $X=I_n$ with $n=1000$. Then, we simulate the response from the linear model \begin{equation*}
  y=X\beta+\epsilon, \epsilon\sim N(0,I_n).
\end{equation*}
The number of relevant features $t$ is set to vary within $\{50,100,200,300,400,500\}$ and the nonzero regression coefficients are equal to $3\sqrt{2\log n}$. 
We set the target FDR level $\alpha = 0.1$ and $\gamma = 0.1$ for F-SLOPE, and set $k=\{5,10,15,20,25,30\}$ and $\alpha = 0.1$ for $k$-SLOPE. Table \ref{tab1} reports the estimation of FDR, $\rm{Prob}\mathrm\{FDP\geq \gamma\}$ and power with 100 repetitions. Figure \ref{orth_kfwer_figure} summaries the results of SLOPE, $k$-SLOPE and F-SLOPE in these trials. These results show our proposed stepdown SLOPEs can reach the FDP control, FDR control and $k$-FWER control flexibly, while SLOPE just can control the FDR. Meanwhile, $k$-SLOPE and F-SLOPE also enjoy the promising power in almost all settings. Furthermore, these experimental results verify the validity of Theorems 4 and 5. Due to the space limitation, we just present the part experimental results (in  Figure \ref{orth_kfwer_figure} and Table \ref{tab1}) and put the comprehensive results in \emph{Supplementary Material C.1}.
\begin{table}[ht]
\centering
\renewcommand\arraystretch{1.15}{
\setlength{\tabcolsep}{2mm}{
\begin{tabular}{c|cccc}
\hline
$t$  & F-SLOPE & $k$-SLOPE & Sd (FDP) & Sd ($k$-FWER) \\ \hline
10 & 0.00    & 0.04    & 0.02       & 0.08            \\
20 & 0.03    & 0.00    & 0.00       & 0.03            \\
30 & 0.00    & 0.01    & 0.01       & 0.01            \\
40 & 0.00    & 0.00    & 0.00       & 0.00            \\
50 & 0.01    & 0.00    & 0.00       & 0.00            \\
60 & 0.00    & 0.00    & 0.00       & 0.00            \\
70 & 0.00    & 0.00    & 0.00       & 0.00            \\
80 & 0.00    & 0.00    & 0.00       & 0.00            \\ \hline
\end{tabular}}}
\caption{$\mathrm{Prob(FDP}>\gamma$) results on the simulated data for multiple mean testing ($k=6$)}\label{tab2}
\end{table}


\subsection{Multiple mean testing from correlated statistics}

\begin{table*}[ht]
\centering
\renewcommand\arraystretch{1.2}{
\setlength{\tabcolsep}{2.1mm}{
\begin{tabular}{c|cccccccccccccccc}
\hline
\multirow{2}{*}{$k$ / $t$} & \multicolumn{8}{c}{weak signals}                              & \multicolumn{8}{c}{moderate signals}                             \\
                     & 10   & 20   & 30   & 40   & 50   & 60   & 70   & 80   & 10   & 20   & 30   & 40   & 50   & 60   & 70   & 80   \\ \hline
2                    & 0.02 & 0.00 & 0.02 & 0.00 & 0.04 & 0.04 & 0.07 & 0.07 & 0.00 & 0.00 & 0.05 & 0.03 & 0.10 & 0.12 & 0.07 & 0.13 \\
4                    & 0.00 & 0.00 & 0.00 & 0.01 & 0.00 & 0.00 & 0.00 & 0.01 & 0.00 & 0.00 & 0.00 & 0.00 & 0.00 & 0.00 & 0.01 & 0.06 \\
6                    & 0.00 & 0.00 & 0.00 & 0.00 & 0.00 & 0.00 & 0.00 & 0.00 & 0.00 & 0.00 & 0.00 & 0.00 & 0.00 & 0.00 & 0.00 & 0.01 \\
8                    & 0.00 & 0.00 & 0.00 & 0.01 & 0.00 & 0.00 & 0.00 & 0.01 & 0.00 & 0.00 & 0.00 & 0.00 & 0.00 & 0.02 & 0.00 & 0.00 \\ \hline
\end{tabular}}}
\caption{$k$-FWER of $k$-SLOPE ($m=2n$) on the simulated data under the weak and moderate signals (different $t$ and $k$).}
\label{tab3}
\end{table*}

\begin{table}[ht]
\centering
\renewcommand\arraystretch{1.1}{
\setlength{\tabcolsep}{4.5mm}{
\begin{tabular}{c|cccc}
\hline
\multirow{2}{*}{$t$} & \multicolumn{2}{c}{$m = 2n$} & \multicolumn{2}{c}{$m = n/2$} \\
                   & weak       & moder       & weak        & moder       \\ \hline
10                 & 0.03       & 0.00        & 0.05        & 0.00        \\
20                 & 0.07       & 0.00        & 0.03        & 0.01        \\
30                 & 0.01       & 0.00        & 0.00        & 0.00        \\
40                 & 0.00       & 0.01        & 0.00        & 0.00        \\
50                 & 0.00       & 0.00        & 0.00        & 0.00        \\
60                 & 0.00       & 0.00        & 0.00        & 0.00        \\
70                 & 0.01       & 0.00        & 0.00        & 0.01        \\
80                 & 0.00       & 0.02        & 0.00        & 0.00        \\ \hline
\end{tabular}}}

\caption{\rm{Prob}(FDP$>\gamma$) of F-SLOPE on the simulated data under the weak and moderate signals (different $t$).}
\label{tab4}
\end{table}
We exemplify the properties of our proposed methods as applied to the typical multiple testing problem with correlated test statistics. Similar to \cite{slope}, we consider the following case. Researchers conduct $n=1000$ experiments in each of $p=5$ randomly selected laboratories. Observation results are modeled as 
\begin{equation*}
y_{i, j}=\mu_{i}+\tau_{j}+z_{i, j}, \quad 1 \leq i \leq n,1 \leq j \leq p,
\end{equation*}
where $\tau_{j}\sim N(0,\sigma_{\tau}^2)$ is the laboratory impact factors, $z_{i, j}\sim N(0,\sigma_{z}^2)$ is the errors and they are independent of each other. Our goal is to test whether $\mu_i$ is equal to 0, i.e. $H_i: \mu_i= 0,i = 1,2,...,n$. Averaging the observed values of 5 laboratories, we get the mean of results 
\begin{equation*}
\bar{y}_{i}=\mu_{i}+\bar{\tau}+\bar{z}_{i}, \quad 1 \leq i \leq n,
\end{equation*}
where $\bar{y}=(\bar{y}_1,...,\bar{y}_n)^T$ is drawn independently from  $N(\mu,\Sigma)$, where $\Sigma_{i, i}=\frac{1}{5}\sigma_{\tau}^2=\rho$ and $\Sigma_{i, j}=\frac{1}{5}(\sigma_{\tau}^2+\sigma_{z}^2)=\sigma^2$ for $i\neq j$ \cite{slope}. 
The key problem is to judge whether the marginal means of a multivariate correlation Gaussian vector disappear or not.
One classical solution is to perform marginal tests with $\bar{y}$ statistic, which depends on the stepdown procedure to control $k$-FWER or FDP \cite{kFWER2}. In other words, we sort the $\bar{y}$ sequence with $|\bar{y}|_{(1)}\geq |\bar{y}|_{(2)}\geq\cdots |\bar{y}|_{(m)}$. Then we use the stepdown procedure  with the $k$-FWER critical values or FDP critical values. Another solution is to ``whiten the noise", i.e., the regression equation is reduced to
\begin{equation}
\tilde{y}=\Sigma^{-1 / 2} \bar{y}=\Sigma^{-1/2} \mu+\epsilon,\label{whiten}
\end{equation}
where $\epsilon\sim N(0,I_p)$, $\Sigma^{-1/2}$ is the regression design matrix. If $\Sigma^{-1/2}$ is closed to the orthogonal matrix, the multiple tests problem is transformed into the feature selection problem under the approximate orthogonal design, where $k$-SLOPE and F-SLOPE can provide better performance.

Similar with \cite{slope}, we set $\sigma_{\tau}^2=\sigma_{z}^2=2.5$ and  consider a sparse setting, where the number of the relevant features $t$ is $\{10,20,30,40,50,60,70,80,90,100\}$. The nonzero mean is set to $2\sqrt{2\log{p}}/c$, where $c$ is equal to the Euclidean norm of each of the columns of $\Sigma^{-1/2}$. We set $\alpha = \gamma = 0.1$ for all FDP controlled methods, and set $k=\{2,4,6,8,10\}$ and $\alpha = 0.1$ for $k$-FWER controlled methods. Figure \ref{figure2} shows the FDR, $k$-FWER and power provided by F-SLOPE, $k$-SLOPE, the stepdown procedures for FDP control (Sd(FDP)), and the stepdown procedures for $k$-FWER control (Sd($k$-FWER)). Table \ref{tab2} shows $\mathrm{Prob(FDP}>\gamma)$ for F-SLOPE, $k$-SLOPE and the stepdown procedures. These experimental results show that our proposed methods ensure the FDP, FDR and $k$-FWER control simultaneously, while Sd (FDP) (or Sd ($k$-FWER)) focuses on  controlling  the FDP (or $k$-FWER) and FDR. However, F-SLOPE and $k$-SLOPE have greater power than the stepdown procedures.
Therefore, our proposed methods have better performance than the classical stepdown procedures in multiple tests. Please refer to \emph{Supplementary Material C.2} for more empirical results.

\subsection{Experiments of Gaussian Design Setting}
\begin{figure}[ht]
	\centering
	\includegraphics[width=1\linewidth]{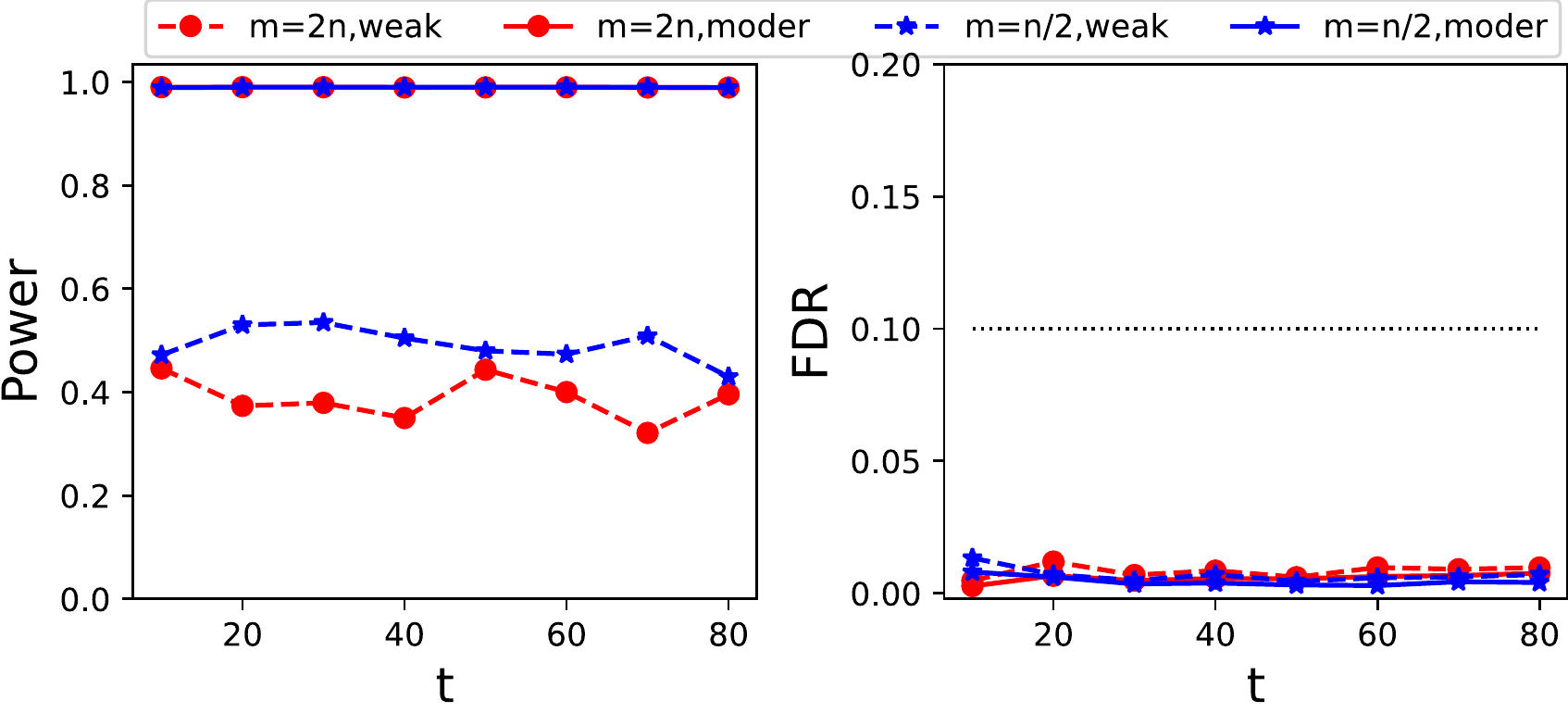}
	\caption{ Power and FDR of F-SLOPE under Gaussian design (different $t$). The black dashed line indicates the target FDR level.}
	\label{figure3}
\end{figure}
We study the performance of $k$-SLOPE and F-SLOPE in general setting. Following the strategy in \cite{slope}, let the entries of the design matrix $X$ are i.i.d $N (0,1/n)$ with $n=5000$. The number of relevant features $t$ varies $\{10,20,30,40,50,60,70,80\}$. Moderate signals having nonzero regression coefficients is set to $2\sqrt{2\log m}$, while this value is set to $\sqrt{2\log m}$ for weak signals. We set $\alpha = \gamma = 0.1$ for F-SLOPE, and set $k=\{2,4,6,8,10\}$ and $\alpha = 0.1$ for $k$-SLOPE.

Then we consider two scenarios: (1) $m=2n$; (2) $m=n/2$. Table \ref{tab4} and Figure \ref{figure3} illustrate F-SLOPE keeps the $\mathrm{Prob(FDP}>\gamma)$ and FDR below the norminal level under both scenarios ($m=2n$ and $m=n/2$), whether the signals are weak and moderate. Meanwhile, Figure \ref{figure3} also shows F-SLOPE ($m=n/2$) has greater power than F-SLOPE ($m=2n$) under weak signals, while F-SLOPE ($m=2n$) and F-SLOPE  ($m=n/2$) have similar power under the moderate signals. As shown in Table \ref{tab3}, $k$-SLOPE control $k$-FWER under both scenarios ($m=2n$ and $m=n/2$). In addition, the power of $k$-SLOPE also has nice performance under the moderate signals. Moreover, experimental results verify the validity of $k$-SLOPE with $\lambda_{k\mathrm{G}^*}$ and F-SLOPE with $\lambda_{\mathrm{FG}^*}$. See \emph{Supplementary Material C.3} for additional experimental results.


\section{Conclusion} \label{sec6}

This paper formulated two  feature selection approaches based on   the SLOPE technique \cite{slope}. Different from the existing works concerning the FDR control, the current models focus on the $k$-FWER control and FDP control for feature selection. With the help of stepdown procedure \cite{kFWER2}, we established their theoretical guarantees under the orthogonal design. Simulated experiments validated the effectiveness of the proposed  stepdown SLOPEs. 

\section{Acknowledgements}
This work was supported in part by National Natural Science Foundation of China under Grant No. 12071166 and by the Fundamental Research Funds for the Central
Universities of China under Grant 2662020LXQD002. We are grateful to the anonymous AAAI reviewers for their constructive comments.

\bibliography{aaai23}

\begin{thebibliography}{43}
\providecommand{\natexlab}[1]{#1}

\bibitem[{Aggarwal and Yadav(2016)}]{aggarwal-2016biofdr}
Aggarwal, S.; and Yadav, A.~K. 2016.
\newblock False discovery rate estimation in proteomics.
\newblock In \emph{Statistical Analysis in Proteomics}, 119--128. Springer.

\bibitem[{Alem{\'a}n et~al.(2017)Alem{\'a}n, Duryea, Guerra, McEwan, Mu{\~n}oz,
  Stampini, and Williamson}]{aleman2017effects}
Alem{\'a}n, X.; Duryea, S.; Guerra, N.~G.; McEwan, P.~J.; Mu{\~n}oz, R.;
  Stampini, M.; and Williamson, A.~A. 2017.
\newblock The effects of musical training on child development: A randomized
  trial of El Sistema in Venezuela.
\newblock \emph{Prevention Science}, 18(7): 865--878.

\bibitem[{Barber and Cand{\`e}s(2015)}]{barber-2015knockoffs}
Barber, R.~F.; and Cand{\`e}s, E.~J. 2015.
\newblock Controlling the false discovery rate via knockoffs.
\newblock \emph{The Annals of Statistics}, 43(5): 2055--2085.

\bibitem[{Barber, Cand{\`e}s, and Samworth(2020)}]{barber-2020knockoffs}
Barber, R.~F.; Cand{\`e}s, E.~J.; and Samworth, R.~J. 2020.
\newblock Robust inference with knockoffs.
\newblock \emph{The Annals of Statistics}, 48(3): 1409--1431.

\bibitem[{Beck and Teboulle(2009)}]{beck2009fast}
Beck, A.; and Teboulle, M. 2009.
\newblock A fast iterative shrinkage-thresholding algorithm for linear inverse
  problems.
\newblock \emph{SIAM journal on imaging sciences}, 2(1): 183--202.

\bibitem[{Becker, Cand{\`e}s, and Grant(2011)}]{becker2011templates}
Becker, S.~R.; Cand{\`e}s, E.~J.; and Grant, M.~C. 2011.
\newblock Templates for convex cone problems with applications to sparse signal
  recovery.
\newblock \emph{Mathematical programming computation}, 3(3): 165--218.

\bibitem[{Bellec, Lecu{\'e}, and Tsybakov(2018)}]{bellec}
Bellec, P.~C.; Lecu{\'e}, G.; and Tsybakov, A.~B. 2018.
\newblock SLOPE meets LASSO: Improved oracle bounds and optimality.
\newblock \emph{The Annals of Statistics}, 46(6B): 3603--3642.

\bibitem[{Benjamini and Hochberg(1995)}]{benjamini-1995bh}
Benjamini, Y.; and Hochberg, Y. 1995.
\newblock Controlling the false discovery rate: a practical and powerful
  approach to multiple testing.
\newblock \emph{Journal of the Royal statistical society: series B
  (Methodological)}, 57(1): 289--300.

\bibitem[{Bogdan et~al.(2015)Bogdan, Van Den~Berg, Sabatti, Su, and
  Cand{\`e}s}]{slope}
Bogdan, M.; Van Den~Berg, E.; Sabatti, C.; Su, W.; and Cand{\`e}s, E.~J. 2015.
\newblock SLOPE—adaptive variable selection via convex optimization.
\newblock \emph{The Annals of Applied Statistics}, 9(3): 1103.

\bibitem[{Breiman(2001)}]{randomforest}
Breiman, L. 2001.
\newblock Random forests.
\newblock \emph{Machine learning}, 45(1): 5--32.

\bibitem[{Brzyski et~al.(2019)Brzyski, Gossmann, Su, and Bogdan}]{Groupslope}
Brzyski, D.; Gossmann, A.; Su, W.; and Bogdan, M. 2019.
\newblock Group slope--adaptive selection of groups of predictors.
\newblock \emph{Journal of the American Statistical Association}, 114(525):
  419--433.

\bibitem[{Brzyski et~al.(2017)Brzyski, Peterson, Sobczyk, Cand{\`e}s, Bogdan,
  and Sabatti}]{brzyski}
Brzyski, D.; Peterson, C.~B.; Sobczyk, P.; Cand{\`e}s, E.~J.; Bogdan, M.; and
  Sabatti, C. 2017.
\newblock Controlling the rate of GWAS false discoveries.
\newblock \emph{Genetics}, 205(1): 61--75.

\bibitem[{Cand{\`e}s et~al.(2018)Cand{\`e}s, Fan, Janson, and
  Lv}]{Cands2016PanningFG}
Cand{\`e}s, E.; Fan, Y.; Janson, L.; and Lv, J. 2018.
\newblock Panning for gold:‘model-X’knockoffs for high dimensional
  controlled variable selection.
\newblock \emph{Journal of the Royal Statistical Society: Series B (Statistical
  Methodology)}, 80(3): 551--577.

\bibitem[{Chen et~al.(2017)Chen, Wang, Deng, and Huang}]{chenh-nips17}
Chen, H.; Wang, X.; Deng, C.; and Huang, H. 2017.
\newblock Group sparse additive machine.
\newblock In \emph{Proceedings of the 31st International Conference on Neural
  Information Processing Systems}.

\bibitem[{Chen et~al.(2021)Chen, Wang, Zheng, Deng, and
  Huang}]{chenh-tnnls2020}
Chen, H.; Wang, Y.; Zheng, F.; Deng, C.; and Huang, H. 2021.
\newblock Sparse modal additive model.
\newblock \emph{IEEE Transactions on Neural Networks and Learning Systems},
  32(6): 2373--2387.

\bibitem[{Delattre and Roquain(2015)}]{delattre2015new}
Delattre, S.; and Roquain, E. 2015.
\newblock New procedures controlling the false discovery proportion via
  Romano--Wolf’s heuristic.
\newblock \emph{The Annals of Statistics}, 43(3): 1141--1177.

\bibitem[{Fan and Lv(2010)}]{fan2010selective}
Fan, J.; and Lv, J. 2010.
\newblock A selective overview of variable selection in high dimensional
  feature space.
\newblock \emph{Statistica Sinica}, 20(1): 101.

\bibitem[{Fan, Demirkaya, and Lv(2019)}]{Fan}
Fan, Y.; Demirkaya, E.; and Lv, J. 2019.
\newblock Nonuniformity of p-values can occur early in diverging dimensions.
\newblock \emph{The Journal of Machine Learning Research}, 20(1): 2849--2881.

\bibitem[{Fan et~al.(2017)Fan, Lyu, Ying, and Hu}]{topk}
Fan, Y.; Lyu, S.; Ying, Y.; and Hu, B.-G. 2017.
\newblock Learning with average top-k loss.
\newblock In \emph{Proceedings of the 31st International Conference on Neural
  Information Processing Systems}, 497--505.

\bibitem[{Farcomeni(2008)}]{Farcomen}
Farcomeni, A. 2008.
\newblock A review of modern multiple hypothesis testing, with particular
  attention to the false discovery proportion.
\newblock \emph{Statistical Methods in Medical Research}, 17(4): 347--388.

\bibitem[{Ferreira and Zwinderman(2006)}]{ferreira-2006bh}
Ferreira, J.; and Zwinderman, A. 2006.
\newblock On the benjamini--hochberg method.
\newblock \emph{The Annals of Statistics}, 34(4): 1827--1849.

\bibitem[{Hammersley and Morton(1954)}]{hammersley1954poor}
Hammersley, J.~M.; and Morton, K.~W. 1954.
\newblock Poor man's monte carlo.
\newblock \emph{Journal of the Royal Statistical Society: Series B
  (Methodological)}, 16(1): 23--38.

\bibitem[{Javanmard and Javadi(2019)}]{javanmard2019false}
Javanmard, A.; and Javadi, H. 2019.
\newblock False discovery rate control via debiased lasso.
\newblock \emph{Electronic Journal of Statistics}, 13(1): 1212--1253.

\bibitem[{Jiang et~al.(2022)Jiang, Bogdan, Josse, Majewski, Miasojedow,
  Ro{\v{c}}kov{\'a}, and Group}]{abslope}
Jiang, W.; Bogdan, M.; Josse, J.; Majewski, S.; Miasojedow, B.;
  Ro{\v{c}}kov{\'a}, V.; and Group, T. 2022.
\newblock Adaptive bayesian SLOPE: Model selection with incomplete data.
\newblock \emph{Journal of Computational and Graphical Statistics}, 31(1):
  113--137.

\bibitem[{Kos and Bogdan(2020)}]{kos}
Kos, M.; and Bogdan, M. 2020.
\newblock On the asymptotic properties of SLOPE.
\newblock \emph{Sankhya A}, 82(2): 499--532.

\bibitem[{Kremer et~al.(2020)Kremer, Lee, Bogdan, and Paterlini}]{kremer}
Kremer, P.~J.; Lee, S.; Bogdan, M.; and Paterlini, S. 2020.
\newblock Sparse portfolio selection via the sorted $l_1$-norm.
\newblock \emph{Journal of Banking \& Finance}, 110: 105687.

\bibitem[{Larsson, Bogdan, and Wallin(2020)}]{strongslope}
Larsson, J.; Bogdan, M.; and Wallin, J. 2020.
\newblock The strong screening rule for SLOPE.
\newblock \emph{Advances in Neural Information Processing Systems}, 33:
  14592--14603.

\bibitem[{Larsson et~al.(2022)Larsson, Klopfenstein, Massias, and
  Wallin}]{larsson}
Larsson, J.; Klopfenstein, Q.; Massias, M.; and Wallin, J. 2022.
\newblock Coordinate descent for SLOPE.
\newblock \emph{arXiv preprint arXiv:2210.14780}.

\bibitem[{Lee, Sobczyk, and Bogdan(2019)}]{lee}
Lee, S.; Sobczyk, P.; and Bogdan, M. 2019.
\newblock Structure learning of Gaussian Markov random fields with false
  discovery rate control.
\newblock \emph{Symmetry}, 11(10): 1311.

\bibitem[{Lehmann and Romano(2005)}]{kFWER2}
Lehmann, E.; and Romano, J.~P. 2005.
\newblock Generalizations of the familywise error rate.
\newblock \emph{The Annals of Statistics}, 1138--1154.

\bibitem[{Lemhadri, Ruan, and Tibshirani(2021)}]{lemhadri-2021lassonet}
Lemhadri, I.; Ruan, F.; and Tibshirani, R. 2021.
\newblock Lassonet: neural networks with feature sparsity.
\newblock In \emph{International Conference on Artificial Intelligence and
  Statistics}, 10--18. PMLR.

\bibitem[{Luo et~al.(2019)Luo, Sun, Toh, and Xiu}]{luo2019solving}
Luo, Z.; Sun, D.; Toh, K.-C.; and Xiu, N. 2019.
\newblock Solving the OSCAR and SLOPE models using a semismooth Newton-based
  augmented Lagrangian method.
\newblock \emph{J. Mach. Learn. Res.}, 20(106): 1--25.

\bibitem[{Nydick(2012)}]{wishart}
Nydick, S.~W. 2012.
\newblock The wishart and inverse wishart distributions.
\newblock \emph{Electronic Journal of Statistics}, 6(1-19).

\bibitem[{Ravikumar et~al.(2009)Ravikumar, Lafferty, Liu, and Wasserman}]{spam}
Ravikumar, P.; Lafferty, J.; Liu, H.; and Wasserman, L. 2009.
\newblock Sparse additive models.
\newblock \emph{Journal of the Royal Statistical Society: Series B (Statistical
  Methodology)}, 71(5): 1009--1030.

\bibitem[{Riccobello et~al.(2022)Riccobello, Bogdan, Bonaccolto, Kremer,
  Paterlini, and Sobczyk}]{riccobello}
Riccobello, R.; Bogdan, M.; Bonaccolto, G.; Kremer, P.~J.; Paterlini, S.; and
  Sobczyk, P. 2022.
\newblock Sparse graphical modelling via the sorted $l_1$-norm.
\newblock \emph{arXiv preprint arXiv:2204.10403}.

\bibitem[{Romano and Shaikh(2006)}]{stepup}
Romano, J.~P.; and Shaikh, A.~M. 2006.
\newblock Stepup procedures for control of generalizations of the familywise
  error rate.
\newblock \emph{The Annals of Statistics}, 34(4): 1850--1873.

\bibitem[{Romano and Wolf(2007)}]{romano2007control}
Romano, J.~P.; and Wolf, M. 2007.
\newblock Control of generalized error rates in multiple testing.
\newblock \emph{The Annals of Statistics}, 35(4): 1378--1408.

\bibitem[{Romano, Sesia, and Cand{\`e}s(2020)}]{romano-2020knockoffs}
Romano, Y.; Sesia, M.; and Cand{\`e}s, E. 2020.
\newblock Deep knockoffs.
\newblock \emph{Journal of the American Statistical Association}, 115(532):
  1861--1872.

\bibitem[{Su and Cand{\`e}s(2016)}]{slopestatics}
Su, W.; and Cand{\`e}s, E. 2016.
\newblock SLOPE is adaptive to unknown sparsity and asymptotically minimax.
\newblock \emph{The Annals of Statistics}, 44(3): 1038--1068.

\bibitem[{Tibshirani(1996)}]{lasso}
Tibshirani, R. 1996.
\newblock Regression shrinkage and selection via the lasso.
\newblock \emph{Journal of the Royal Statistical Society: Series B
  (Methodological)}, 58(1): 267--288.

\bibitem[{Van~der Laan, Dudoit, and Pollard(2004)}]{vanderLaanDudoitPollard}
Van~der Laan, M.~J.; Dudoit, S.; and Pollard, K.~S. 2004.
\newblock Augmentation procedures for control of the generalized family-wise
  error rate and tail probabilities for the proportion of false positives.
\newblock \emph{Statistical Applications in Genetics and Molecular Biology},
  3(1).

\bibitem[{Yu, Kaufmann, and Lederer(2021)}]{yu-2021biofdr}
Yu, L.; Kaufmann, T.; and Lederer, J. 2021.
\newblock False discovery rates in biological networks.
\newblock In \emph{International Conference on Artificial Intelligence and
  Statistics}, 163--171. PMLR.

\bibitem[{Zhao et~al.(2022)Zhao, Chen, Wang, Li, Gong, Wang, and
  Zheng}]{zhao-2022error}
Zhao, X.; Chen, H.; Wang, Y.; Li, W.; Gong, T.; Wang, Y.; and Zheng, F. 2022.
\newblock Error-based knockoffs inference for controlled feature selection.
\newblock In \emph{Proceedings of the AAAI Conference on Artificial
  Intelligence}, volume~36, 9190--9198.

\end{thebibliography}

\clearpage
\section{Supplementary Material A}
To improve the readability, we summarize the main notations of this paper in Table \ref{notations}.

\begin{table*}[ht]
\centering
\renewcommand\arraystretch{1}{
\setlength{\tabcolsep}{12mm}{
\begin{tabular}{cl}
\toprule
{Notations} & Descriptions\\
\midrule
$n$ & the sample size\\
$m$ & the dimension of input\\
$\mathcal{X}, \mathcal{Y}$ & the input space $\mathcal{X}\subset\mathbb{R}^m$ and the output space $\mathcal{Y}\subset \mathbb{R}$, respectively \\
$X$ & the design matrix with $X:=(X_1,X_2...,X_n)^T\subset \mathbb{R}^{n\times m}$\\
$y$ & the observation vector with $y=(y_1,y_2,...,y_n)^T\in \mathbb{R}^n$\\
$\beta$ & the coefficient vector with $\beta=(\beta_1,\beta_2,...,\beta_n)\in \mathbb{R}^n$\\
$\epsilon$ & the Gaussian error with zero mean and the known variance $\sigma^2I_n$\\
$V$ & the number of false selected features\\
$R$ & the number of identified features\\
$\lambda$ & the m-dimensional regularization parameter vector with $\lambda_1\geq\lambda_2\geq...\geq\lambda_m\geq0$\\
$m_0$ & the number of true null hypotheses\\
$\alpha$ & the desired FDR, FDP and $k$-FWER level\\
$\gamma$ & the fixed level which FDP exceeds \\
$k$ & the at most number of false selected features tolerated \\
\bottomrule
\end{tabular}}}
\caption{Notations}
\label{notations}
\end{table*}

\section{Supplementary Material B}
Before providing the proofs of Theorems 4 and 5, we first establish stepping stones including Lemma 1-3. The following lemmas (See the proof of Lemma 1 and 2 in \cite{slope} and the proof of Lemma 3 in \cite{topk}) is key to these theorems.
\begin{lemma}\cite{slope}\label{lemma1}
	Let $H_i$ be a null hypotheses and let $r\geq 1$. Then 
	\begin{equation*}
		\left\{y: H_{i} \text { is rejected and } R=r\right\}=\left\{y:\left|y_{i}\right|>\lambda_{r} \text { and } R=r\right\}.
	\end{equation*}
\end{lemma}

\begin{remark}
Lemma 1 demonstrates that when the number of selected features $R$ is equal to $r$ ($r\geq1$), the conditions for $H_i$ to be rejected are its corresponding $|y_i|$ is greater than $\lambda_r$.
\end{remark}

\begin{lemma}\cite{slope}\label{lemma2}
	Consider applying the SLOPE procedure to $\tilde{y}=(y_1,...,y_{i-1},y_{i+1},...y_p)$ with weight $\tilde{\lambda}=(\lambda_2,...,\lambda_p)$ and let $\tilde{R}$ be the number of rejections this procedure makes. Then with $r>1$,
\begin{equation*}
    \{y:|y_i|>\lambda_r\text{ and } R=r\}\subset\{y:|y_i|>\lambda_r\text{ and }\tilde{R}=r-1\}.
\end{equation*}
\end{lemma}

\begin{remark}
Lemma 2 explains that if $|y_i|$ is arbitrarily removed, the number of features selected is larger than the previous set.
\end{remark}

\begin{lemma}\cite{topk}\label{lemma3}
	Let $s_{[k]}$ be the top-k element of a set $S=\{s_1,...,s_n\}$, such as $s_{[1]}\geq s_{[2]}\geq...\geq s_{[n]}$. $\sum_{i=1}^ks_{[i]}$ is a convex function of $(s_1,...,s_n)$. Furthermore, for $s_i\geq0$ and $i=1,...,n$, we have $\sum_{i=1}^ks_{[i]}=\min_{\lambda\geq0}\{k\lambda+\sum_{i=1}^n[s_i-\lambda]_+\}$, where $[a]_+=\max\{0,a\}$ is the hinge function.
\end{lemma}

\begin{remark}
Lemma 3 gives the equivalent form of the sum of the first $k$ non-negative values. 
\end{remark}

\subsection{Proof of Theorem 4}
Suppose that $X$ is the orthogonal design matrix and $\sigma^2 = 1$. Now we have $\tilde{y}=X'y\sim N(\beta,I_m)$ and $m$ null hypotheses $H_i:\beta_i=0, i=1,...,m$, in which the first $m_0$ hypotheses are true null hypotheses. Then order the absolute value of $y$ corresponding to the first $m_0$ true null hypotheses; denote them 
\begin{equation*}
    |\hat{q}|_{(1)}\geq|\hat{q}|_{(2)}\geq...\geq|\hat{q}|_{(m_0)}.
\end{equation*}
Let $k$ be at most the number of false selected features tolerated or the number of indices that satisfy $\beta_i\neq 0$ in $\{1,...,m_0\}$. Assume $m_0\geq k$ or there is nothing to prove. We have,
\begin{equation}
    k\text{-}\mathrm{FWER}=\mathrm{Prob}(V\geq k)=\sum_{r=1}^m \mathrm{Prob}(V\geq k \text{ and } R=r).\label{eq1}
\end{equation}
Through the stepdown procedure \cite{kFWER2} and Lemma \ref{lemma1} \cite{slope}, $k$-SLOPE commits at least $k$ false rejections if and only if 
\begin{equation*}
	|\hat{q}|_{(1)}\geq \lambda_r,|\hat{q}|_{(2)}\geq\lambda_r,...,|\hat{q}|_{(k)}\geq\lambda_r,
\end{equation*}
when the number of selected features $R$ is equal to $r$. Then
\begin{align*}
    &\mathrm{Prob}(V\geq k \text{ and } R=r)\\
    &=\mathrm{Prob}(|\hat{q}|_{(1)}\geq\lambda_r,...,|\hat{q}|_{(k)}\geq\lambda_r \text{ and } R=r)\\
    &\leq \mathrm{Prob}(|\hat{q}|_{(k)}\geq\lambda_r \text{ and } R=r).
\end{align*}
Due to the Lemma 2 \cite{slope} and the independence between $|\hat{q}|_{(k)}$ and $\tilde{y}$, we have
\begin{align}
    &\mathrm{Prob}(V\geq k \text{ and } R=r)\\&\leq \mathrm{Prob}(|\hat{q}|_{(k)}\geq\lambda_r)\mathrm{Prob}(\tilde{R}=r-1)\label{eq2}
\end{align}
It is not difficult to find $\mathrm{Prob}(|\hat{q}|_{(i)}\geq\lambda_r)$ non-increasing, so 
\begin{align}
  \mathrm{Prob}(|\hat{q}|_{(k)}\geq\lambda_r) \leq \frac{1}{k}\sum_{i=1}^k\mathrm{Prob}(|\hat{q}|_{(i)}\geq\lambda_r)\label{eq3}
\end{align}
Next, combined with lemma \ref{lemma3} \cite{topk}, we have that
\begin{align*}
    &\frac{1}{k}\sum_{i=1}^k\mathrm{Prob}(|\hat{q}|_{(i)}\geq\lambda_r)\\&=\frac{1}{k}\min_{t\geq0}\{kt+\sum_{i=1}^{m_0}[\mathrm{Prob}(|\hat{q}_i|\geq\lambda_r)-t]_+\}\\
    &\leq\min_{t_0\leq t\leq \alpha}\{t+\frac{m_0}{k}[t_0-t]_+\},
\end{align*}
where $t_0=\mathrm{Prob}(|\hat{q}_i|\geq\lambda_r)=k\alpha/(p+k-r)$. Plugging these inequalities into (\ref{eq1}) gives
\begin{align*}
    k\text{-}\mathrm{FWER} &= \sum_{r=1}^m\mathrm{Prob}(|\hat{q}|_{(1)}\geq\lambda_r,...,|\hat{q}|_{(k)}\geq\lambda\text{ and } R=r)\\
    &\leq\sum_{r=1}^m\min_{t_0\leq t\leq \alpha}\{t+\frac{p_0}{k}[t_0-t]_+\}\mathrm{Prob}(\tilde{R}=r-1)\\
    &\leq \sum_{r=1}^m \frac{k\alpha}{p+k-r} \mathrm{Prob}(\tilde{R}=r-1)\leq \alpha,
\end{align*}
which completes the proof of Theorem 4.

\subsection{Proof of Theorem 5}
It is the same as the assumption of Theorem 4. Let the number of true hypotheses is non-zero, i.e. $m_0>0$, otherwise it is not necessary to prove it again. Given $\gamma\in(0,1)$, we have
\begin{align}
    \mathrm{Prob}(\mathrm{FDP}>\gamma)&=\sum_{r=1}^m\mathrm{Prob}(\frac{V}{R}>\gamma\text{ and } R=r)\\
    &=\sum_{r=1}^m\mathrm{Prob}(V\geq k(R) \text{ and }R=r),\label{eq4}
\end{align}
where $k(R)=\lfloor{\gamma R}\rfloor+1$. We observed that (\ref{eq2}) is similar to (\ref{eq1}) and the only difference is whether the value of $k$ is affected by the number of selected features. Based on Lemma \ref{lemma1} \cite{slope} and the stepdown procedure \cite{kFWER2}, 
\begin{align*}
    &\mathrm{Prob}(V\geq k(R) \text{ and } R=r)\\&=\mathrm{Prob}(|\hat{q}|_{(1)}\geq\lambda_r,...,|\hat{q}|_{(k(R))}\geq\lambda_r\text{ and }R=r)\\
    &\leq \mathrm{Prob}(|\hat{q}|_{(k)}\geq\lambda_r \text{ and } R=r).
\end{align*}
Analogy to (\ref{eq2}) and (\ref{eq3}), we have
\begin{align*}
    &\mathrm{Prob}(V\geq k(R) \text{ and } R=r)\leq\mathrm{Prob}(|\hat{q}|_{(k(R))}\geq\lambda_r)\mathrm{Prob}(\tilde{R}=r)\\
    &\leq\frac{1}{k(R)}\sum_{i=1}^{k(R)}\mathrm{Prob}(|\hat{q}|_{(i)}\geq\lambda_r)\mathrm{Prob}(\tilde{R}=r-1)
\end{align*}
Then combined with lemma \ref{lemma3} \cite{topk}, we have 
\begin{align*}
     &\frac{1}{k(R)}\sum_{i=1}^{k(R)}\mathrm{Prob}(|\hat{q}|_{(i)}\geq\lambda_r)\\&=\frac{1}{k(R)}\min_{t\geq0}\{k(R)\cdot t+\sum_{i=1}^{m_0}[\mathrm{Prob}(|\hat{q}_i|\geq\lambda_r)-t]_+\}\\
    &\leq\min_{t_0\leq t\leq \alpha}\{t+\frac{p_0}{k(R)}[t_0-t]_+\},
\end{align*}
where $t_0=(k(r)+1)q/(p+k(r)+1-r)$. Plugging these inequalities into (\ref{eq2}) gives
\begin{align*}
    &\mathrm{Prob}(FDP>\gamma)\\&=\sum_{r=1}^m{P}(V\geq k(R)\text{ and } R=r)\\
    &\leq\sum_{r=1}^m \min_{t_0\leq t\leq q}\{t+\frac{m_0}{k(R)}[t_0-t]_+\}\mathrm{Prob}(\tilde{R}=r-1)\\
    &\leq \sum_{r=1}^m\frac{k(r)q}{p+k(r)-r}\mathrm{Prob}(\tilde{R}=r-1)\leq \alpha,
\end{align*}
which finishes the proof of Theorem 5.

\clearpage

\section{Supplementary Material C}

\subsection{C.1 Orthogonal design}
Table \ref{tab11} shows $k$-SLOPE keeps the FDR and $\mathrm{Prob(FDP>}\gamma)$ below the standard level and has good power under the orthogonal design.

\subsection{C.2 Multiple mean testing from correlated statistics}

Table \ref{tab12} and \ref{tab13} illustrate $k$-SLOPE and the stepdown procedure ($k$-FWER) \cite{kFWER2} ensure the FDR and FDP control, while the power of $k$-SLOPE is greater than that of the stepdown procedure ($k$-FWER) \cite{kFWER2} in multiple testing.

\subsection{C.3 Gaussian design}
\begin{figure}[ht]
	\centering
	\includegraphics[width=1\linewidth]{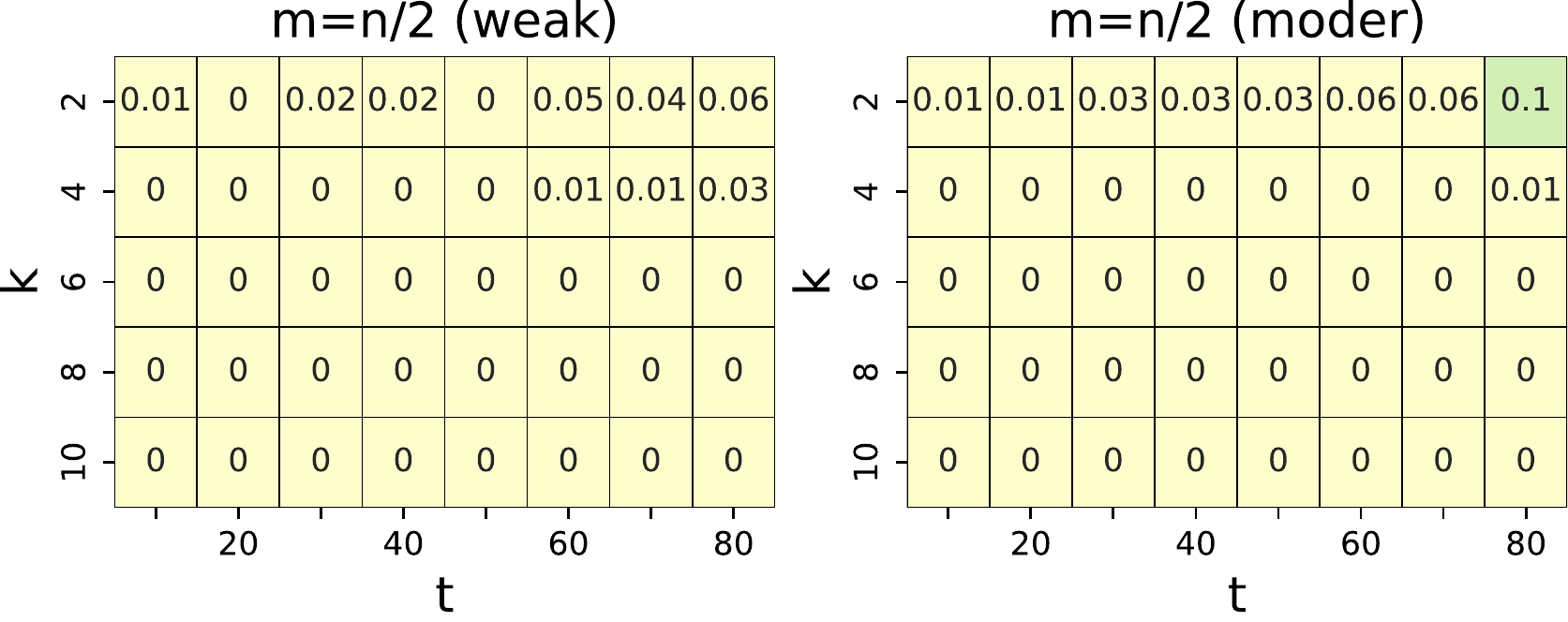}
	\caption{$k$-FWER for $k$-SLOPE ($m=n/2$) under the weak and moderate signals (different t and $k$).}
	\label{fig4}
\end{figure}

Figure \ref{fig4} and \ref{fig6} and Table \ref{tab15} and \ref{tab16} show important properties of $k$-SLOPE under the Gaussian design. Figure 2 illustrates $k$-FWER provided by F-SLOPE under the weak and moderate signals. These results explain $k$-SLOPE has nice performance for the FDR and $k$-FWER control. The power of $k$-SLOPE under the moderate signals is greater than that under the weak signals. 

\begin{figure}[ht]
	\centering
	\includegraphics[width=1\linewidth]{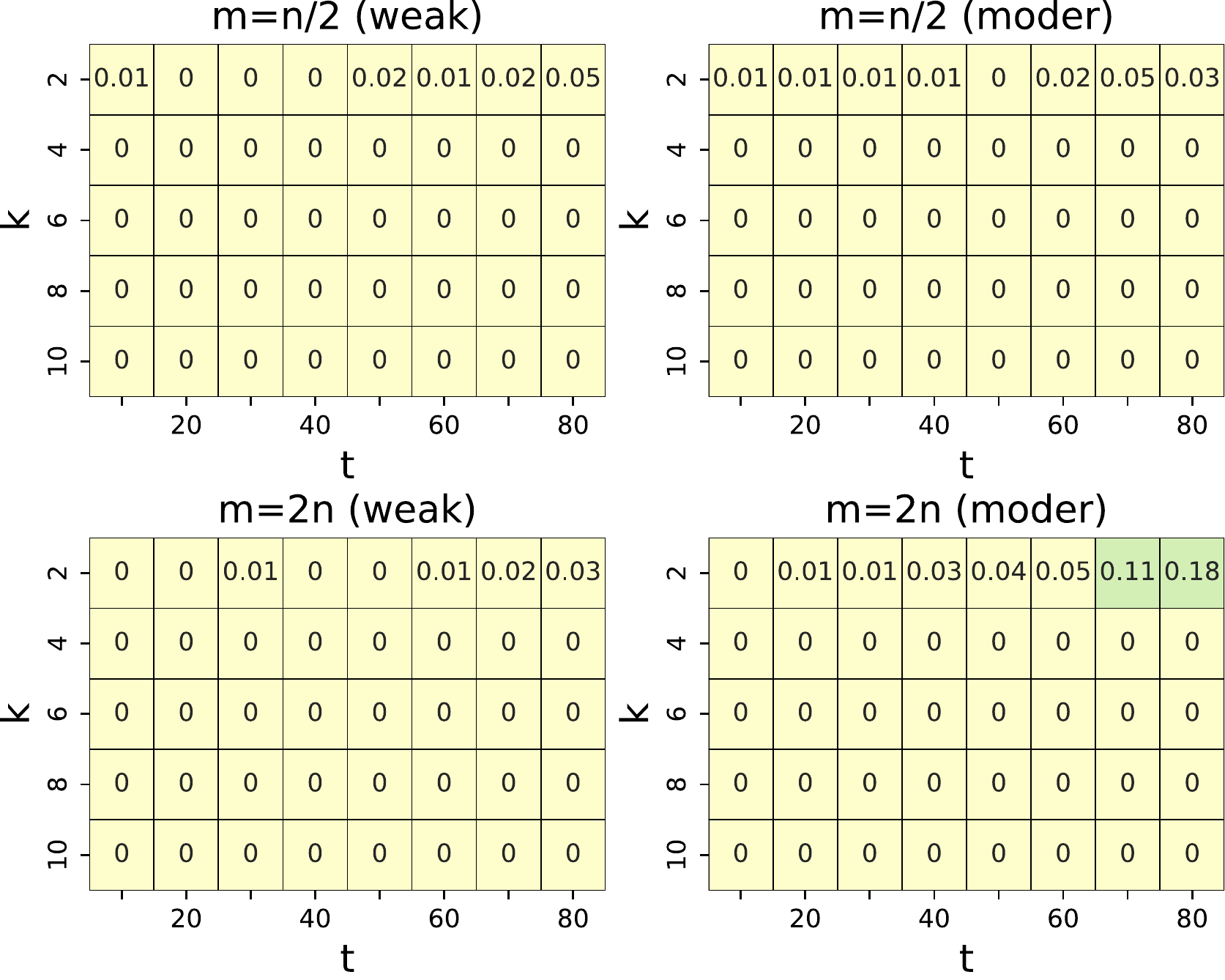}
	\caption{$k$-FWER for F-SLOPE under the weak and moderate signals  (different t and $k$).}
	\label{fig5}
\end{figure}
\begin{figure}[htb]
	\centering
	\includegraphics[width=1\linewidth]{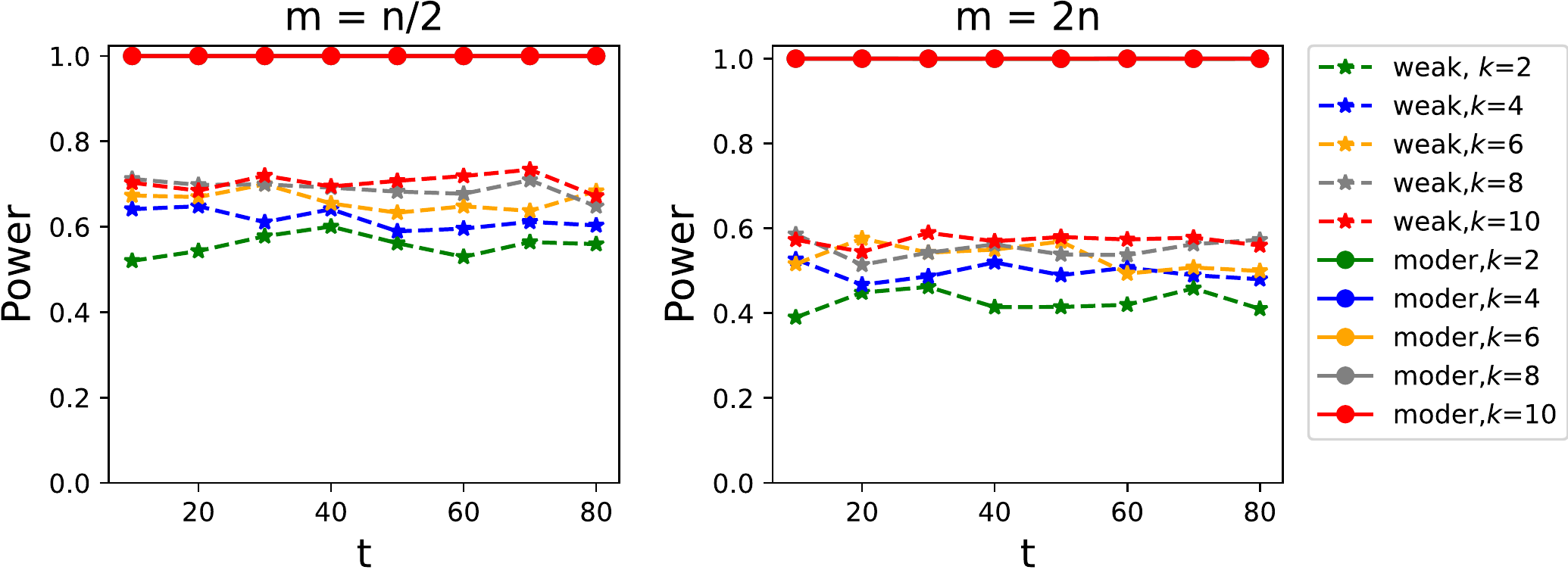}
	\caption{Power of $k$-SLOPE on the simulated data (different t and $k$).}
	\label{fig6}
\end{figure}

\begin{table*}[ht]
\centering
\renewcommand\arraystretch{1.2}{
\setlength{\tabcolsep}{0.85mm}{
\begin{tabular}{c|cccccc|cccccc|cccccc}
\hline
\multirow{2}{*}{t / $k$} & \multicolumn{6}{c|}{$\mathrm{Prob(FDP>}\gamma)$}    & \multicolumn{6}{c|}{FDR}& \multicolumn{6}{c}{Power}\\
                   & 5    & 10   & 15   & 20   & 25   & 30   & 5    & 15   & 20   & 25   & 30   & 5   & 10    & 15   & 20   & 25   & 400   & 500   \\ \hline
50                 & 0.001 & 0.000 & 0.001 & 0.000 & 0.000 & 0.010 & 0.006 & 0.010 & 0.015 & 0.019 & 0.026 & 0.032 & 1.000 & 1.000 & 1.000 & 1.000 & 0.007 & 1.000 \\
100                & 0.000 & 0.000 & 0.000 & 0.002 & 0.000 & 0.000 & 0.002 & 0.005 & 0.007 & 0.010 & 0.013 & 0.014 & 0.998 & 1.000 & 0.990 & 1.000 & 1.000 & 1.000 \\
200                & 0.001 & 0.001 & 0.000 & 0.000 & 0.000 & 0.002 & 0.001 & 0.003 & 0.004 & 0.005 & 0.006 & 0.007 & 1.000 & 1.000 & 1.000 & 1.000 & 1.000 & 0.998 \\
300                & 0.002 & 0.000 & 0.000 & 0.000 & 0.003 & 0.000 & 0.001 & 0.001 & 0.002 & 0.003 & 0.004 & 0.005 & 1.000 & 1.000 & 1.000 & 0.996 & 1.000 & 1.000 \\
400                & 0.000 & 0.000 & 0.001 & 0.000 & 0.000 & 0.000 & 0.001 & 0.001 & 0.002 & 0.002 & 0.003 & 0.004 & 0.995 & 1.000 & 0.994 & 1.000 & 0.005 & 1.000 \\
500                & 0.000 & 0.000 & 0.000 & 0.000 & 0.000 & 0.000 & 0.001 & 0.001 & 0.002 & 0.002 & 0.002 & 0.003 & 0.997 & 1.000 & 1.000 & 1.000 & 1.000 & 1.000 \\ \hline
\end{tabular}}}
\caption{Results of $k$-SLOPE under the orthogonal design (different t and $k$.)}
\label{tab11}
\end{table*}
\begin{table*}[htbp]
\centering
\renewcommand\arraystretch{1.1}{
\begin{tabular}{c|ccccc|ccccc|ccccc}
\hline
\multirow{2}{*}{t / k} & \multicolumn{5}{c|}{FDR}              & \multicolumn{5}{c|}{$\mathrm{Prob(FDP>}\gamma)$} & \multicolumn{5}{c}{Power}             \\
                     & 2     & 4     & 6     & 8     & 10    & 2      & 4      & 6      & 8      & 10     & 2     & 4     & 6     & 8     & 10    \\ \hline
10                   & 0.011 & 0.022 & 0.022 & 0.026 & 0.043 & 0.00   & 0.04   & 0.04   & 0.05   & 0.07   & 0.927 & 0.956 & 0.961 & 0.974 & 0.966 \\
20                   & 0.004 & 0.003 & 0.010 & 0.013 & 0.020 & 0.00   & 0.00   & 0.00   & 0.00   & 0.01   & 0.940 & 0.948 & 0.963 & 0.968 & 0.976 \\
30                   & 0.003 & 0.005 & 0.010 & 0.010 & 0.013 & 0.00   & 0.00   & 0.01   & 0.00   & 0.00   & 0.928 & 0.944 & 0.960 & 0.965 & 0.973 \\
40                   & 0.001 & 0.004 & 0.004 & 0.007 & 0.007 & 0.00   & 0.00   & 0.00   & 0.00   & 0.00   & 0.910 & 0.941 & 0.955 & 0.956 & 0.964 \\
50                   & 0.000 & 0.001 & 0.003 & 0.003 & 0.005 & 0.00   & 0.00   & 0.00   & 0.00   & 0.00   & 0.910 & 0.942 & 0.953 & 0.954 & 0.970 \\
60                   & 0.001 & 0.002 & 0.003 & 0.003 & 0.003 & 0.00   & 0.00   & 0.00   & 0.00   & 0.00   & 0.910 & 0.928 & 0.945 & 0.960 & 0.963 \\
70                   & 0.001 & 0.001 & 0.002 & 0.002 & 0.004 & 0.00   & 0.00   & 0.00   & 0.00   & 0.00   & 0.900 & 0.927 & 0.946 & 0.953 & 0.959 \\
80                   & 0.000 & 0.001 & 0.001 & 0.002 & 0.001 & 0.00   & 0.00   & 0.00   & 0.00   & 0.00   & 0.890 & 0.924 & 0.939 & 0.945 & 0.955 \\ \hline
\end{tabular}}
\caption{Results of $k$-SLOPE on the simulated data (different t and $k$).}
\label{tab12}
\end{table*}
\begin{table*}[htbp]
\renewcommand\arraystretch{1}{
\setlength{\tabcolsep}{2.1mm}{
\begin{tabular}{c|ccccc|ccccc|ccccc}
\hline
\multirow{2}{*}{t / $k$} & \multicolumn{5}{c|}{FDR}              & \multicolumn{5}{c|}{$\mathrm{Prob(FDP>}\gamma)$} & \multicolumn{5}{c}{Power}             \\
                     & 2     & 4     & 6     & 8     & 10    & 2      & 4      & 6      & 8      & 10     & 2     & 4     & 6     & 8     & 10    \\ \hline
10                   & 0.012 & 0.038 & 0.047 & 0.075 & 0.077 & 0.04   & 0.11   & 0.08   & 0.02   & 0.24   & 0.916 & 0.949 & 0.963 & 0.966 & 0.967 \\
20                   & 0.015 & 0.022 & 0.035 & 0.040 & 0.048 & 0.01   & 0.03   & 0.03   & 0.07   & 0.07   & 0.928 & 0.942 & 0.950 & 0.963 & 0.970 \\
30                   & 0.008 & 0.021 & 0.022 & 0.024 & 0.036 & 0.00   & 0.00   & 0.01   & 0.01   & 0.02   & 0.913 & 0.944 & 0.949 & 0.955 & 0.958 \\
40                   & 0.007 & 0.011 & 0.021 & 0.028 & 0.032 & 0.00   & 0.00   & 0.00   & 0.02   & 0.03   & 0.906 & 0.942 & 0.951 & 0.957 & 0.961 \\
50                   & 0.006 & 0.014 & 0.018 & 0.025 & 0.027 & 0.00   & 0.00   & 0.00   & 0.00   & 0.03   & 0.909 & 0.935 & 0.937 & 0.954 & 0.954 \\
60                   & 0.007 & 0.010 & 0.018 & 0.019 & 0.027 & 0.00   & 0.00   & 0.00   & 0.00   & 0.00   & 0.899 & 0.919 & 0.935 & 0.944 & 0.955 \\
70                   & 0.006 & 0.014 & 0.018 & 0.024 & 0.025 & 0.00   & 0.00   & 0.00   & 0.00   & 0.00   & 0.879 & 0.916 & 0.928 & 0.942 & 0.951 \\
80                   & 0.006 & 0.010 & 0.016 & 0.020 & 0.024 & 0.00   & 0.00   & 0.00   & 0.00   & 0.00   & 0.868 & 0.900 & 0.923 & 0.940 & 0.943 \\ \hline
\end{tabular}}}
\caption{Results of the stepdown procedure ($k$-FWER) on the simulated data (different t and $k$).}
\label{tab13}
\end{table*}


\begin{table*}[ht]
\centering
\begin{tabular}{c|cccccccc}
\hline
\multirow{3}{*}{k / t} & \multicolumn{4}{c}{weak}                                  & \multicolumn{4}{c}{moderate}                              \\
                   & 20           & 40           & 60           & 80           & 20           & 40           & 60           & 80           \\ \hline
2                  & 0.016, 0.014 & 0.009, 0.012 & 0.011, 0.014 & 0.007, 0.014 & 0.010, 0.005  & 0.006, 0.008 & 0.006, 0.009 & 0.006, 0.009 \\
4                  & 0.015, 0.027 & 0.017, 0.023 & 0.013, 0.021 & 0.014, 0.020  & 0.018, 0.015 & 0.009, 0.013 & 0.009, 0.013 & 0.009, 0.017 \\
6                  & 0.030, 0.035  & 0.02, 0.024  & 0.018, 0.029 & 0.014, 0.023 & 0.018, 0.017 & 0.015, 0.022 & 0.010, 0.021  & 0.009, 0.021 \\
8                  & 0.039, 0.054 & 0.022, 0.032 & 0.021, 0.032 & 0.019, 0.033 & 0.027, 0.026 & 0.023, 0.023 & 0.015, 0.019 & 0.013, 0.023 \\
10                 & 0.043, 0.06  & 0.029, 0.036 & 0.023, 0.037 & 0.022, 0.045 & 0.032, 0.034 & 0.024, 0.023 & 0.015, 0.025 & 0.018, 0.028 \\ \hline
\end{tabular}
\caption{FDR of $k$-SLOPE under the weak and moderate signals (different t and k).}
\label{tab15}
\end{table*}

\begin{table*}[ht]
\centering
\renewcommand\arraystretch{1.1}{
\setlength{\tabcolsep}{3.6mm}{
\begin{tabular}{c|cccccccc}
\hline
\multirow{2}{*}{k / t} & \multicolumn{4}{c}{weak}                       & \multicolumn{4}{c}{moderate}          \\
                   & 20         & 40         & 60         & 80      & 20         & 40      & 60      & 80   \\ \hline
2                  & 0.04, 0.02 & 0.00, 0.00       & 0.00, 0.00       & 0.00, 0.00    & 0.00, 0.00       & 0.00, 0.00    & 0.00, 0.00    & 0.00, 0.00 \\
4                  & 0.02, 0.09 & 0.00, 0.08    & 0.01, 0.00    & 0.00, 0.00    & 0.02, 0.00    & 0.00, 0.00    & 0.00, 0.00    & 0.00, 0.00 \\
6                  & 0.05, 0.10  & 0.00, 0.03    & 0.00, 0.02    & 0.00, 0.00    & 0.01, 0.00    & 0.00, 0.00    & 0.00, 0.00    & 0.00, 0.00 \\
8                  & 0.07, 0.24 & 0.03, 0.03 & 0.02, 0.02 & 0.00, 0.01 & 0.01, 0.00    & 0.01, 0.00 & 0.00, 0.01 & 0.00, 0.00 \\
10                 & 0.15, 0.23 & 0,02, 0.06 & 0, 0.02    & 0.00, 0.07 & 0.05, 0.05 & 0.00, 0.00    & 0.00, 0.01 & 0.00, 0.00 \\ \hline
\end{tabular}}}
\caption{$\mathrm{Prob(FDP>}\gamma)$ of $k$-SLOPE under the weak and moderate signals (different t and k)}.
\label{tab16}
\end{table*}

\end{document}